\documentclass[12pt,twoside,fleqn,reqno]{amsart}
\usepackage{amssymb,palatino,graphics,hyperref}%
\usepackage[mathcal]{euler}
\usepackage[T1]{fontenc}

\makeatletter
\renewcommand{\section}{\@startsection {section}{1}{\z@}%
                                   {-3.5ex \@plus -1ex \@minus -.2ex}%
                                   {.5\linespacing}%
                                   {\normalfont\scshape\centering}}
\@ifpackageloaded{euler}{
 \def\B@R#1#2{\raisebox{-.07ex}{$#1#2$}\mkern-6mu}
 \renewcommand{\hbar}{{\mspace{1mu}\mathpalette\B@R{\mathchar'26}h}}
  \renewcommand{\cong}{\stackrel{\raise1pt\hbox{$\sim$}}{\boldsymbol{\smash=}}}
  }{}
\makeatother

\numberwithin{equation}{section}

\newtheorem{thm}{Theorem}[section]
\newtheorem{lem}[thm]{Lemma}
\newtheorem{cor}[thm]{Corollary}
\newtheorem{prop}[thm]{Proposition}
  
\theoremstyle{definition}
\newtheorem{definition}{Definition}[section]

\def\beq#1\eeq{\begin{equation}#1\end{equation}}
\newcommand{\A}{\mathcal{A}}
\renewcommand{\k}{\hbar}
\newcommand{\Ak}{\A_\k}
\renewcommand{\AA}{\mathbb{A}}
\newcommand{\C}{\mathcal{C}}
\newcommand{\Ci}{\C^\infty}
\newcommand{\M}{\mathcal{M}}
\newcommand{\Mt}{\widetilde\M}
\newcommand{\co}{\mathbb{C}}
\newcommand{\R}{\mathbb{R}}
\DeclareMathOperator{\Image}{Im}
 \renewcommand{\Im}{\Image}
\newcommand{\inner}{\mathbin{\raise0.3ex\hbox{$\lrcorner$}}}
\newcommand{\Li}{\mathcal{L}}
\newcommand{\D}{\mathcal{D}}
\newcommand{\Hi}{\mathcal{H}}
\newcommand{\cs}{\mbox{\upshape C}\ensuremath{{}^*}} 
\newcommand{\abs}[1]{\lvert#1\rvert}
\newcommand{\Abs}[1]{\left|#1\right|}
\newcommand{\norm}[1]{\lVert#1\rVert}

\renewcommand{\S}{\mathcal{S}}
\newcommand{\T}{\mathbb{T}}
\newcommand{\G}{\Gamma}
\DeclareMathOperator{\SU}{SU}
\DeclareMathOperator{\End}{End}
\newcommand{\g}{\mathfrak{g}}
\newcommand{\Or}{\mathcal{O}}

\newcommand{\isom}{\mathrel{\widetilde\longrightarrow}}
\newcommand{\Z}{\mathbb{Z}}
\DeclareMathOperator{\tr}{tr}
\newcommand{\Om}{\mbox{\boldmath$\Omega$}}
\newcommand{\Po}{\mathcal{P}}
\newcommand{\F}{\mathcal F}
\newcommand{\Fp}{\F^\perp}
\newcommand{\Wedge}{\wedge}
\newcommand{\Dabrowski}{D¡browski} 
\DeclareMathOperator{\rk}{rk}
\newcommand{\lcon}{\mbox{\boldmath$\nabla$}}
\DeclareMathOperator{\Ad}{Ad}
\newcommand{\normal}{\mathrel{\triangleleft}}
\DeclareMathOperator{\Isom}{Isom}
\newcommand{\ta}{\mathfrak t}
\newcommand{\Lie}{\mathcal{L}}
\newcommand{\B}{\mathcal{B}}

\newcommand{\Sn}{S}
\newcommand{\Bk}{\B_\k}
\newcommand{\Nil}{\mathrm{Nil}^3}
\newcommand{\Poincare}{Poincar\'{e}}
\newcommand{\Podles}{Podle\'{s}}
\newcommand{\Mreg}{\M_{\mathrm{reg}}}
\newcommand{\op}{\mathrm{op}}
\DeclareMathOperator{\Tr}{Tr}
\newcommand{\BB}{\mathbb{B}}
\newcommand{\iso}{\mathop{\mathfrak{iso}}}
\newcommand{\so}{\mathop{\mathfrak{so}}}
\newcommand{\tnabla}{\widetilde\nabla}
\newcommand{\Birkhauser}{Birkh\"auser}
 
\title{The Structure of Noncommutative Deformations}
\author{Eli Hawkins}
\subjclass[2000]{58B34; \emph{Secondary} 46L65, 53D17}
\keywords{Noncommutative Geometry, Deformation Quantization, Poisson Geometry}

\begin{document}
\maketitle
\begin{center}
\vspace{-4ex}
\emph{\small Department of Mathematics}\\
\emph{\small The University of Western Ontario}\\
\emph{\small London, Ontario N6A\,5B7, Canada}\\
{\small mrmuon@mac.com}\\
\end{center}
\begin{abstract}
Noncommutatively deformed geometries, such as the noncommutative torus, do not exist generically. I showed in a previous paper that the existence of such a deformation implies compatibility conditions between the classical metric and the Poisson bivector (which characterizes the noncommutativity). Here I present another necessary condition: the vanishing of a certain rank $5$ tensor. In the case of a compact Riemannian manifold, I use these conditions to prove that the Poisson bivector can be constructed locally from commuting Killing vectors. 
\end{abstract}

\section{Introduction}
The notion of a noncommutative deformation is not uncommon in contemporary mathematics. Quantum mechanics can be viewed as a noncommutative deformation of classical mechanics. Quantum groups are noncommutative deformations of Lie groups. The \Podles\ spheres and other quantum group homogeneous spaces are deformations. The noncommutative torus is a deformation of a torus. $E$-theory is a bivariant $K$-theory constructed from equivalence classes of deformations. Strict and formal deformation quantization are two mathematical settings for studying deformations.

If the starting point of a noncommutative deformation is a geometrical space, then it is natural to try to view the deformation geometrically.  At best, we can try to promote a 1-parameter family of algebras to a 1-parameter family of noncommutative geometries.
In \cite{haw}, I showed that this is not always possible. There are obstructions to deforming integration, $1$-forms, or a metric.  

Any noncommutative deformation of a smooth manifold is characterized by a Poisson structure (geometrically, an antisymmetric bivector field). These obstructions are expressed in the language of Poisson geometry. In particular, the necessary condition for deformation of $1$-forms is the existence of a flat, torsion-free ``contravariant connection''. The first main result (Sec.~\ref{meta}) in this paper is an obstruction to deforming higher degree differential forms. This obstruction is a tensor $M^{ijk}_{lm}$ which I call the ``metacurvature'' of the flat, torsion-free contravariant connection.

The most general fully fledged example of a noncommutative geometric deformation was given by Connes and Landi \cite{c-l}. Their construction applies to a compact Riemannian manifold with a torus group acting by isometries. An invariant Poisson structure on the torus determines the deformation. This raises the question of how much more general is the most general deformation.

Using all the available obstructions, I show in Section \ref{local} that a noncommutative deformation of a compact Riemannian manifold corresponds to a Poisson structure that can be expressed locally in terms of commuting Killing vectors. Globally, the Poisson structure is induced by an invariant Poisson structure on a group of isometries of a covering space. Using this structure, I sketch a generalization of the Connes-Landi construction. This construction indicates that my obstructions are not only necessary but \emph{sufficient} conditions for the existence of a noncommutative geometric deformation.

One of my obstructions is a condition for deforming integration into a trace on the noncommutative algebras. This is independent of the other obstructions, so it is possible to consider the other obstructions alone. Classifying the solutions to this weaker set of conditions is more difficult in general, but I investigate the simplest case of 2 dimensions and find that the only solutions are a flat torus with a constant Poisson bivector and a round sphere with the Poisson structure associated to the \Podles\ ``standard'' sphere. The latter case corresponds to a noncommutative geometric deformation constructed by \Dabrowski\ and Sitarz; this example satisfies a weakened version of Connes' axioms for noncommutative geometry.

In the remainder of this section, I present further background. In Section \ref{meta}, I derive the metacurvature tensor as the obstruction to the existence of a differential graded Poisson algebra; I prove a formula for the metacurvature in the simplest case (a symplectic manifold). In Section \ref{2-D}, I use this formula to classify noncommutative deformations of compact 2-dimensional Riemannian manifolds, temporarily disregarding integration. In Section \ref{divergence.sec}, I reformulate the obstruction regarding integration. In Section \ref{realizations}, I present several lemmata showing how the metacurvature and related structures behave in a symplectic realization. These are tools for the proofs in Section \ref{Riemannian} where I show that the Poisson bivector for a deformation of a compact Riemannian manifold can be decomposed locally into products of Killing vectors. Globally, I construct a symmetry group from such a Poisson structure and show how any such Poisson structure can be constructed. After some more background, in Section \ref{converse} I sketch the construction of a geometric deformation for any compatible Poisson structure on a compact Riemannian manifold. Finally, in Section \ref{examples}, I give a few simple examples of such Poisson structures on Riemannian manifolds.

\subsection{Notation}
$\C^\infty(\M)$ will denote the space of smooth (infinitely differentiable), $\co$-valued functions on a smooth manifold $\M$. $\G(\M,V)$ will denote the space of smooth sections of a vector bundle $V$ over $\M$. $\Omega^p(\M) := \G(\M,\Wedge^p T^*\M)$ will denote the space of smooth differential $p$-forms, and $\Omega^p(\M,V)$ the space of smooth $p$-forms with coefficients in $V$.

When discussing deformations directly, it is necessary to use complex functions and sections. However, the connections and brackets that I work with all preserve real sections.

I will mostly use index-free notation for tensors. However, it is occasionally necessary to resort to index notation. In index notation, a vertical bar denotes a covariant derivative. Summation is implicit over any repeated index.

If $\F$ is a foliation of $\M$, then $T\F$ denotes the tangent bundle to the foliation; this is a bundle over $\M$. The differential forms along $\F$ are $\Omega^*(\F) := \G(\M,\Wedge^*T^*\F)$.

A vector field acts on functions as a directional derivative operator as in $X(f)$. Multivectors (vectors, bivectors, \emph{et cetera}) are sections of the exterior powers $\Wedge^*T\M$ of the tangent bundle. I will use the symbol $\inner$ to denote not only the contraction of a vector into a differential form, but also the contraction of a multivector into a differential form. This is such that, for instance, $(X\wedge Y)\inner\epsilon = Y\inner(X\inner\epsilon) = \epsilon(X,Y,\dots)$. An exponent on a form or multivector always denotes an exterior power. 

\subsection{Lie Algebroids}
Several of the principal structures that I use here are unified by the concepts of a Lie algebroid and a connection with respect to a Lie algebroid.

\begin{definition}
A \emph{Lie algebroid} is a vector bundle $V$ over a smooth manifold $\M$ with a vector bundle homomorphism $\iota : V \to T\M$ (the \emph{anchor map}) and a Lie algebra structure on $\G(\M,V)$ such that:
\begin{enumerate}
\item
For $v,w\in\G(\M,V)$ and $f\in\Ci(\M)$ 
\[
[v,fw] = f\, [v,w] + \iota_v(f) w
\mbox;
\]
\item
The anchor map intertwines the $V$-bracket with the Lie bracket of vector fields
\[
[\iota_v,\iota_w] = \iota_{[v,w]}
\mbox.
\]
\end{enumerate}
\end{definition}

The tangent bundle is itself a Lie algebroid with the identity as anchor map. A Lie algebra is precisely a Lie algebroid over a point. The holomorphic tangent bundle of a complex manifold is a Lie algebroid. The tangent bundle to a foliation is a Lie algebroid; the anchor being the inclusion map into $T\M$.

The most important Lie algebroid here is the cotangent bundle to a Poisson manifold. In this case, I denote the anchor as $\#:T^*\M\to T\M$; it is defined by
\[
(\#\alpha)^i = \pi^{ji}\alpha_j
\mbox.
\]
The bracket is the Koszul bracket which I denote as $[\;\cdot\;,\;\cdot\;]_\pi$; it is uniquely defined by the defining properties of a Lie algebroid and the condition that
\[
[df,dg]_\pi = d\{f,g\}
\mbox.\]

The usual notion of a connection for a vector bundle naturally generalizes to Lie algebroids. 
\begin{definition}
A \emph{connection} on a vector bundle $W$ with respect to a Lie algebroid $V$ is a map $\nabla$ from sections of $V$ to first order differential operators on $\G(\M,W)$ such that for any $f\in\Ci(\M)$, $v\in\G(\M,V)$, and $w\in\G(\M,W)$:
\[
\nabla_{fv}w = f\,\nabla_vw
\]
and
\[
\nabla_v(fw) = f\,\nabla_vw + \iota_v(f) w
\mbox.\]
\end{definition}

The definition of curvature for such a connection is formally identical to the usual definition:
\[
K(u,v) := \nabla_u\nabla_v - \nabla_v\nabla_u - \nabla_{[u,v]}
\mbox.
\]
It is simple to check that the definitions imply this is a tensor; specifically, it is a section of $\Wedge^2V^*\otimes\End(W)$.

In particular, if $W=V$, then we can define torsion as
\[
T(v,w) := \nabla_vw - \nabla_wv - [v,w]
\mbox.
\]
This is also a tensor, a section of $\Wedge^2V^*\otimes V$. The connections of this kind appearing in this paper are all torsion-free. 

\begin{definition}
Given a foliation $\F$, a \emph{partial connection} on $\F$ is a connection on $T\F$ with respect to $T\F$.
Given a Poisson manifold, a \emph{contravariant connection} \cite{fer1,vai1} is a connection with respect to $T^*\M$; I denote a contravariant connection as $\D$.
\end{definition}

\subsection{Deformations and Previous Results}
Let $\A_0$ be an algebra. 
\begin{definition}
\label{deformation}
A \emph{deformation} of $\A_0$ is an algebraic extension of the form,
\[
0 \to \k\AA \to \AA \stackrel{\Po}\to \A_0 \to 0
\]
where $\k$ is a central multiplier of $\AA$, and for $a\in\AA$, 
\[
\k^2 a =0 \implies a\in\k\AA
\mbox.
\]
\end{definition}
This definition is the weakest possible one for the purposes of this paper; it is essentially equivalent to the definition I gave in \cite{haw}. One could easily replace the last assumption with the simpler but marginally more restrictive condition that $\k$ is not a zero-divisor:
\[
\k a =0 \implies a=0
\mbox.
\]

If $\A_0=\Ci(\M)$, then we should think of $\AA$ as the algebra of smooth functions on a larger noncommutative space. This is a sort of noncommutative cobordism.

If $\AA$ is a deformation of $\Ci(\M)$, then this definition allows us to extract a Poisson bracket from the commutator in $\AA$. This must satisfy the Jacobi identity and is given geometrically by a Poisson bivector, $\pi\in\G(\M,\Wedge^2 T\M)$. 

I now summarize my principal results from \cite{haw}. 

Let $\epsilon\in\Omega^n(\M)$ be a volume form. If integration by $\epsilon$ on $\M$ can be smoothly deformed to a trace, then $\pi$ and $\epsilon$ must satisfy the compatibility condition
\beq
\label{divergence}
d(\pi\inner\epsilon)=0
\mbox.
\eeq

If $\Omega^1(\M)$ and the gradient map $d:\Ci(\M)\to\Omega^1(\M)$ can be smoothly deformed, then there exists a flat, torsion-free contravariant connection $\D$ on $\M$.

If $\M$ has a Riemannian metric and this is extended into a deformed real spectral triple, then the above contravariant connection must be compatible with the metric in the sense that the contravariant derivative of the metric is $0$. This condition can be motivated in other ways, but it requires a sufficiently well defined notion of noncommutative geometry. In fact I will give an improved derivation if this condition in Section \ref{Riemannian}.

For a given metric and Poisson structure there exists a unique torsion-free contravariant connection compatible with the metric. I call this the metric contravariant connection. In this way, the last two conditions can be restated as: The metric contravariant connection is flat.

In the next section I will motivate and define one more condition. Given a flat, torsion-free contravariant connection, there exists a rank 5 tensor which I call the \emph{metacurvature}. If there exists a deformation of differential forms in all degrees, then the metacurvature must be $0$.

\section{Metacurvature}
\label{meta}
Suppose that $0 \to \k\Om^* \to \Om^* \stackrel{\Po}\to \Omega^*(\M) \to 0$
is a deformation of $\Omega^*(\M)$, the differential graded algebra of differential forms. Because $\Omega^*(\M)$ is (graded) commutative, the graded commutator vanishes,
\[
0 = [\sigma,\rho] := \sigma\wedge\rho - (-)^{\deg\sigma\deg\rho} \rho\wedge\sigma
\]
for $\sigma,\rho\in\Omega^*(\M)$. For $\hat\sigma,\hat\rho\in\Om^*$, this implies that $[\hat\sigma,\hat\rho] \in \k \Om^*$. Because of this we can define a generalized Poisson bracket on $\Omega^*(\M)$ by
\[
 \{\Po(\hat\sigma),\Po(\hat\rho)\} = \Po \left(\tfrac1{i\k} [\hat \sigma,\hat \rho]\right)
\mbox.
\]

\begin{thm}
This defines a bracket on $\Omega^*(\M)$, making it a differential graded Poisson algebra.
\end{thm}
\begin{definition}
\label{DG.Poisson}
A \emph{differential graded Poisson algebra} is a graded vector space $\Omega^*$ with 3 operations: $d$, $\wedge$, and $\{\,\cdot\,,\,\cdot\,\}$, such that:
\begin{enumerate}
\item
$d : \Omega^k \to \Omega^{k+1}$ is linear and a differential: $d^2=0$.
\item
$(\Omega^*,\wedge)$ is an associative, graded-commutative algebra.
\item
$(\Omega^*,\{\,\cdot\,,\,\cdot\,\})$ is a graded Lie algebra, i.e., $\{\,\cdot\,,\,\cdot\,\}$ is bilinear, degree $0$, antisymmetric
\beq
\{\sigma,\rho\} = -(-)^{\deg\sigma\deg\rho} \{\rho,\sigma\}
\mbox,
\eeq
and satisfies the graded Jacobi identity,
\beq
\label{Jacobi}
\{\{\sigma,\rho\},\lambda\} = \{\sigma,\{\rho,\lambda\}\} - (-)^{\deg\sigma\deg\rho} \{\rho,\{\sigma,\lambda\}\}
\mbox.
\eeq
\item
$(\Omega^*,d,\wedge)$ is a differential graded algebra, i.e., it also satisfies the Leibniz identity,
\beq
\label{Leibniz}
d(\sigma\wedge\rho) = d\sigma \wedge \rho + (-)^{\deg\sigma} \sigma \wedge d\rho
\mbox.
\eeq
\item
$(\Omega^*,d,\{\,\cdot\,,\,\cdot\,\})$ is a differential graded Lie algebra, i.e., it also satisfies the Leibniz identity
\beq
\label{Leibniz2}
d\{\sigma,\rho\} = \{d\sigma,\rho\} + (-)^{\deg\sigma} \{\sigma,d\rho\}
\mbox.
\eeq
\item
$(\Omega^*,\wedge,\{\,\cdot\,,\,\cdot\,\})$ is a (graded) Poisson algebra, i.e., it also satisfies the product identity\footnote{The terms ``Leibniz identity'' and ``product identity'' are usually synonymous. However, I am using them differently here in order to make a lexical distinction between \eqref{Leibniz2} and \eqref{product}.}
\beq
\label{product}
\{\sigma,\rho\wedge\lambda\} = \{\sigma,\rho\}\wedge\lambda + (-)^{\deg\sigma\deg\rho} \rho\wedge \{\sigma,\lambda\}
\mbox.
\eeq
\end{enumerate}
\end{definition}
\begin{proof}
Of course, $\Omega^*(\M)$ has the structure of a differential graded commutative algebra. We need to check the claims pertaining to the Poisson bracket. It is quite straightforward to check that $\Om^*$ with the commutator satisfies Def.~\ref{DG.Poisson} (except for graded commutativity). The Poisson bracket will inherit these identities from the commutator, provided that it is well defined.

Let $\hat\sigma,\hat\rho\in\Om^*$, $\sigma = \Po(\hat\sigma)$, and $\rho = \Po(\hat\rho)$; we need to check that $\{\sigma,\rho\} := \Po\left(\tfrac1{i\k} [\hat\sigma,\hat\rho]\right)$ is uniquely determined by $\sigma$ and $\rho$. The only way to change $\hat\sigma$ without changing $\sigma$ is to add something of the form $\k \hat\lambda$. However,
\[
[\hat\sigma + \k\hat\lambda,\hat\rho] = [\hat\sigma,\hat\rho] + \k [\hat\lambda,\hat\rho]
\]
which is only changed by an element of $\k^2\Om^*$. Therefore the bracket is well-defined in $\sigma$, and by symmetry, in $\rho$.
\end{proof}

In degree $0$, these structures make $\Ci(\M)$ a Poisson algebra. Hence $\M$ is in particular a Poisson manifold.

This generalized Poisson bracket should not be confused with the Koszul bracket. The Koszul bracket can be naturally extended to differential forms, but it is quite different from a Poisson bracket. The Koszul bracket is of degree $-1$ and is uniquely determined by $\pi$; the Poisson bracket is of degree $0$ and is not determined by $\pi$. For example, if $\alpha,\beta\in\Omega^1(\M)$, then $[a,\beta]_\pi = -[\beta,\alpha]_\pi \in \Omega^1(\M)$ but $\{\alpha,\beta\}=\{\beta,\alpha\}\in \Omega^2(\M)$.

Differential forms with the Koszul bracket form a type of differential Gerstenhaber algebra (the terminology varies); this satisfies identities very similar to those of a differential graded Poisson algebra, but the degree shift in the bracket changes some signs. Grabowski \cite{gra} has constructed a bracket which is determined by $\pi$ and extends the ordinary Poisson bracket; it is a graded Lie bracket (without any degree shift) but does not satisfy the Leibniz or product identities. On the other hand, several authors (e.g., \cite{m-v}) have considered brackets which only fail to satisfy the Leibniz identity \eqref{Leibniz2}.

The product identities \eqref{product} for two functions and a $1$-form imply that the Poisson bracket of a function and $1$-form is given by a contravariant connection $\D$ on $T^*\M$ as
\[
\{f,\alpha\} = \D_{df}\alpha
\mbox.
\]
The Jacobi identity \eqref{Jacobi} for two functions and a $1$-form imply that $\D$ is flat:
\[
\begin{split}
K(df,dg)\alpha &= \D_{df}\D_{dg}\alpha - \D_{dg}\D_{df}\alpha - \D_{[df,dg]_\pi} \alpha \\
&= \D_{df}\D_{dg}\alpha - \D_{dg}\D_{df}\alpha - \D_{d\{f,g\}} \alpha \\
&= \{f,\{g,\alpha\}\} - \{g,\{f,\alpha\}\} - \{\{f,g\},\alpha\} 
= 0
\mbox.
\end{split}
\]
The Leibniz identity \eqref{Leibniz2} for two functions implies that $\D$ is torsion-free:
\[
\begin{split}
T(df,dg) &= \D_{df}dg - \D_{dg}df - [df,dg]_\pi \\
&= \D_{df}dg - \D_{dg}df - d\{f,g\} 
=0
\mbox.
\end{split}
\]

The contravariant connection on $T^*\M$ naturally extends to the exterior powers, $\wedge^kT^*\M$. This is compatible with the exterior product in the obvious way, and so the product identity \eqref{product} implies that the Poisson bracket of a function and differential form is given by $\D$ as
\[
\{f,\sigma\} = \D_{df}\sigma
\mbox.
\]

Now consider the identities involving the bracket of two $1$-forms. The product identity \eqref{product} for $f\in\Ci(\M)$ and $\alpha,\beta\in\Omega^1(\M)$ reads
\[
\{f\alpha,\beta\} = f\,\{\alpha,\beta\} + \alpha\wedge \{f,\beta\}
\mbox.
\]
The Leibniz identity \eqref{Leibniz2} is in this case
\[
d\{f,\alpha\} = \{df,\alpha\} + \{f,d\alpha\}
\mbox.\]
These identities uniquely determine the Poisson bracket of two $1$-forms. However, there is still one identity to be satisfied: the Jacobi identity,
\beq
\begin{split}
0 &\stackrel?= \{f,\{\beta,\gamma\}\} - \{\{f,\beta\},\gamma\} - \{\{f,\gamma\},\beta\} \\
&\;= \D_{df}\{\beta,\gamma\} - \{\D_{df}\beta,\gamma\} - \{\D_{df}\gamma,\beta\}
\mbox.
\end{split}
\label{Jacobi1}
\eeq

Without assuming that this is satisfied, we can consider the properties of the right hand side of eq.~\eqref{Jacobi1}. 
\begin{thm}
A flat, torsion-free contravariant connection determines a tensor $M^{ijk}_{lm}$ symmetric in the contravariant indices and antisymmetric in the covariant indices, such that
\beq
M(df,\beta,\gamma) = \{f,\{\beta,\gamma\}\} - \{\{f,\beta\},\gamma\} - \{\{f,\gamma\},\beta\}
\label{Mdef}
\eeq
if $M$ is viewed as a trilinear map from $1$-forms to $2$-forms.
\end{thm}
\begin{definition}
$M$ is the \emph{metacurvature}.
\end{definition}
\begin{proof}
Begin by taking eq.~\eqref{Mdef} as the definition of some trilinear map, and note that the right hand side of eq.~\eqref{Mdef} is explicitly symmetric in $\beta$ and $\gamma$. We need to check that it is $\Ci(\M)$-linear in either of these arguments. 
For any $g\in\Ci(\M)$,
\[
\begin{split}
M(df,g\beta,\gamma)-g\,M(df,\beta,\gamma)
&= \{f,g\}\,\{\beta,\gamma\} - \{\{f,g\}\beta,\gamma\} \\
&\qquad + \{f,\beta\wedge\{g,\gamma\}\} - \{f,\beta\}\wedge\{g,\gamma\} \\
&\qquad - \beta\wedge\{g,\{f,\gamma\}\} \\
&= \beta\wedge \left(-\{\{f,g\},\gamma\} 
+  \{f,\{g,\gamma\}\} 
- \{g,\{f,\gamma\}\}\right) \\
&= 0
\mbox.
\end{split}
\]
The first two steps use the product identities; the last step uses the lower degree Jacobi identity (the flatness of $\D$).

Now, consider the Jacobi identity 
\[
0 = \{\{f,g\},dh\} - \{f,\{g,dh\}\} + \{g,\{f,dh\}\}
\mbox.\]
Applying $d$ to this equation and using the Leibniz identities for $d$ shows that 
\[
M(df,dg,dh) = \{f,\{dg,dh\}\} - \{\{f,dg\},dh\} - \{\{f,dh\},dg\}
\]
is symmetric under the exchange of $f$ and $g$. Since it is $\Ci(\M)$-linear in $dg$ it must also be $\Ci(\M)$-linear in $df$. This shows that the right hand side of eq.~\eqref{Mdef} is indeed given by a $\Ci(\M)$-trilinear map from $1$-forms to $2$-forms. Such a map is equivalent to a tensor $M^{ijk}_{lm}$ such that $M(\alpha,\beta,\gamma)_{lm} = M^{ijk}_{lm}\alpha_i\beta_j\gamma_k$.
\end{proof} 

\begin{thm}
\label{same}
The following are equivalent:
\begin{enumerate}
\item
A generalized Poisson bracket making $\Omega^*(\M)$ a differential graded Poisson algebra.
\item
A Poisson structure on $\M$ and a flat, torsion-free contravariant connection with $M=0$.
\end{enumerate}
\end{thm}
\begin{proof}
We have already seen that the first structure determines the second. 

The product identity, $\{f\sigma,\rho\} = f\,\{\sigma,\rho\} + \sigma \wedge \D_{df}\rho$, shows that a generalized Poisson bracket $\{\sigma,\rho\}$ is first order differential in both arguments, so it can be constructed in a coordinate chart.
Let $x^i$ be the coordinates. Decomposing a differential form in these coordinates simply means writing it as a sum of products of functions and the basis $1$-forms $dx^i$. So, we can compute any bracket using the product identity \eqref{product} and the fundamental brackets:
\begin{gather*}
\{f,g\}\,\mbox,\\
\{f,dx^i\} = \D_{df}dx^i \mbox,\\
\intertext{and}
\{dx^i,dx^j\} = d\left(\D_{dx^i}dx^j\right) \mbox.
\end{gather*}
This bracket satisfies all the product identities \eqref{product} because these are consistent with associativity. It satisfies the Leibniz identities \eqref{Leibniz2} because they are consistent with the product identities.

This leaves the Jacobi identities \eqref{Jacobi}. Using the coordinate decomposition and the product identities, any Jacobi identity reduces down to the Jacobi identities involving functions and exact $1$-forms. These are satisfied because $\pi$ is Poisson, $\D$ is flat, $M=0$, and
\begin{align*}
\{\{df,dg\},dh\} - \{df,\{dg,dh\}\} - \{dg,\{df,dh\}\}
&= - d\left[M(df,dg,dh)\right] \\
&= 0
\end{align*}
\end{proof}

If the Poisson bivector $\pi$ is invertible, then its inverse is a symplectic $2$-form. 
In that case the flat, torsion-free, contravariant connection $\D$ is related to a flat, torsion-free, covariant connection $\nabla$ by $\#\D_\alpha\beta = \nabla_{\#\alpha}\#\beta$.
\begin{thm}
\label{symplectic.M}
If $\pi=\omega^{-1}$ then,
\begin{subequations}
\begin{align}
M^{ijk}_{lm} &= \D^i\D^j\D^k\omega_{lm} \label{contravariant.M}\\
&= -\pi^{ai}\pi^{bj}\pi^{ck}\omega_{dl}\omega_{em}  \pi^{de}_{\;|abc} \label{covariant.M}
\end{align}
\end{subequations}
\end{thm}
\begin{proof}
Consider any point $x\in\M$.
Because $\nabla$ is flat, any vector in $T_x\M$ extends to a $\nabla$-constant vector field in a neighborhood of $x$. Equivalently, any covector in $T^*_x\M$ extends to a $\D$-constant $1$-form in a neighborhood of $x$. Because of this and symmetry, in order to compute $M$ at $x$, it is sufficient to compute $M(\alpha,\alpha,\alpha)$ for any $1$-form $\alpha\in\Omega^1(U)$ with $\D\alpha=0$ defined over some neighborhood $U\ni x$.

Because $\D\alpha=0$, we have $\{f,\alpha\}=0$ for any $f\in\C^\infty(\M)$. Using this, we compute
\[
\begin{split}
\{f\,dg,\alpha\} &= f\,\{dg,\alpha\} + dg\wedge\{f,\alpha\} = f\,\{dg,\alpha\}\\
&= f\left(\{dg,\alpha\} - d\{g,\alpha\}\right) = - f\,\{g,d\alpha\}\\
&= -f\, \D_{dg} d\alpha
= - \D_{f dg} d\alpha
\end{split}
\]
which shows that $\{\beta,\alpha\} = - \D_\beta d\alpha$ for any $\beta\in\Omega^1(U)$.

Only the first term of \eqref{Mdef} survives in
\[
\begin{split}
M(df,\alpha,\alpha) &= \{f,\{\alpha,\alpha\}\}
= \D_{df} \{\alpha,\alpha\}
= -\D_{df} \D_\alpha d\alpha
\end{split}
\]
Because $M$ is just a tensor, this gives
\[
M(\alpha,\alpha,\alpha) = - \D_\alpha^2 d\alpha .
\]

The vector field $\#\alpha$ is covariantly constant ($\nabla\#\alpha=0$) so the Lie derivative is equal to the covariant derivative, $\Lie_{\#\alpha}=\nabla_{\#\alpha}$. Applying this to $\omega$ gives,
\[
\nabla_{\#\alpha}\omega = \Lie_{\#\alpha}\omega = d(\#\alpha\inner\omega) = d\alpha
\mbox.
\]
The map $\#$ naturally extends to differential forms as in $\#\omega = -\pi$. Because $\pi$ is inverse to $\omega$, its derivative can be rewritten as $\nabla_{\#\alpha}\omega = \#^{-1}\nabla_{\#\alpha}\pi$. 

This gives the expressions,
\begin{align}
\label{M.symplectic}
M(\alpha,\alpha,\alpha) &= - \D_\alpha^2 \#^{-1}\nabla_{\#\alpha}\pi \\
&= -\D_\alpha^3\#^{-1}\pi =  \D_\alpha^3 \omega .
\end{align}
This gives eq.~\eqref{contravariant.M}. To rewrite eq.~\eqref{M.symplectic} in terms of $\nabla$, consider
\[
\# M(\alpha,\alpha,\alpha) = - \#\D_\alpha^2 \#^{-1}\nabla_{\#\alpha}\pi 
= - \nabla_{\#\alpha}^3 \pi .
\]
This gives eq.~\eqref{covariant.M}. 
\end{proof}

This means that for a symplectic manifold, $M=0$ if and only if $\pi$ is quadratic in the locally affine structure defined by the flat connection.

Equation \eqref{contravariant.M} suggests an analogue of the Bianchi identity for the metacurvature. In the symplectic case, the metacurvature is the third contravariant derivative of the symplectic $2$-form. In general, it behaves \emph{as if} it is the third derivative of a $2$-form.
\begin{prop}
$\D^i M^{jkl}_{mn}$ is totally symmetric in the contravariant indices.
\end{prop}
\begin{proof}
It is sufficient to prove that $\D^i M^{jkl}_{mn} =\D^j M^{ikl}_{mn}$.

We can rewrite eq.~\eqref{Mdef} slightly as
\[
M(df,\beta,\gamma) = \D_{df}\{\beta,\gamma\} - \{\D_{df}\beta,\gamma\} - \{\beta,\D_{df}\gamma\}
\mbox.
\]
Formally, this expresses $M$ as the contravariant derivative of the Poisson bracket on two $1$-forms. The proof of the symmetry of the second derivative is formally the same as if $\{\,\cdot\,,\,\cdot\,\}$ were a $\Ci(\M)$-bilinear map. 
\end{proof}

\section{2 Dimensions}
\label{2-D}
Suppose that a Riemannian manifold is deformed into a noncommutative geometry. In particular, assume that differential forms are deformed in a way that is compatible with the metric. By the results of \cite{haw} and Thm.~\ref{same}, this means that the metric contravariant connection is flat and has vanishing metacurvature, $M=0$.

\Dabrowski\ and Sitarz \cite{d-s} have constructed an interesting example of noncommutative geometry on the \Podles\ ``standard'' sphere, a noncommutative deformation of $S^2$. This is particularly interesting because it satisfies some, but not all, of Connes' axioms for noncommutative geometry. Because of this, the homological and spectral dimensions are not $2$ and there is no trace corresponding to integration on $S^2$.

A spectral triple satisfying all of Connes axioms has a canonical trace that plays the role of integration. If a Riemannian manifold is smoothly deformed into such spectral triples, then the volume form and Poisson structure must satisfy the compatibility condition \eqref{divergence}. However, this condition is quite independent of the conditions for the deformation of differential forms. This independence, and the \Dabrowski-Sitarz example, suggest that it would be interesting to set aside eq.~\eqref{divergence}. This makes the analysis much more difficult in general, but it is tractable in 2 dimensions.
 
\begin{thm}
Let $\M$ be a compact, connected, $2$-dimensional Riemannian manifold with a nonzero Poisson structure. Suppose that there exists a deformation of differential forms compatible with the metric. Then $\M$ is either
a torus with a constant metric and Poisson bivector
or
a sphere with constant curvature and the Poisson structure corresponding to the \Podles\ standard sphere.  
\end{thm}
\begin{proof}
We need to solve the conditions that the metric contravariant connection be flat and $M=0$. By assumption, the Poisson bivector $\pi$ is not identically $0$, therefore there is (at least) an open submanifold $\Sigma := \{x\in\M \mid \pi(x)\neq 0\}$. 

Where $\pi$ does not vanish, it is invertible, so $\Sigma$ is symplectic. Over $\Sigma$, $\#$ intertwines $\D$ with a covariant connection $\nabla$ and the metric $ds^2$ with another metric $ds'^2$; the metric tensors are related by
\[
g^{ij} = \pi^{ik}\pi^{jl}g'_{kl}
\mbox.
\]
The facts that $\D$ is flat, torsion free, and compatible with $ds^2$ imply that $\nabla$ is the Levi-Civita connection of $ds'^2$ and is flat.

A bivector in 2 dimensions has only one independent component, so we can write $\pi$ in the form
\[
\pi = h {\epsilon'}^{-1}
\]
where $h$ is a scalar function and $\epsilon'\in\Omega^2(\Sigma)$ is the volume form of $ds'^2$. The relationship between the metrics is therefore a conformal rescaling, $g_{ij} = h^{-2} g'_{ij}$.

By direct computation,
\[
\pi \inner \epsilon = -h^{-1}
\]
is (the restriction of) a continuous function on a compact manifold; therefore it is bounded and there is a constant $A$ such that $0<A\leq \abs{h}$ over $\Sigma$.

Because $(\Sigma,ds'^2)$ is flat, we can cover it with Cartesian coordinate charts. 
By Thm.~\ref{symplectic.M}, $M=0$ implies that $\pi$ is quadratic in any Cartesian chart. Since the volume form $\epsilon'$ is constant,  this just means that $h$ is a quadratic function in any Cartesian chart.

Let $u:[0,L)\to\Sigma$ be a $ds'^2$ geodesic of length $L$. The function $h\circ u(t)$ is quadratic and because $t$ is bounded, there exists another constant $B$ such that,
\[
0 < A \leq \abs{h\circ u(t)} \leq B .
\]
The length of $u$ in $ds^2$ is finite:
\[
\int_0^L \frac{dt}{\abs{h\circ u(t)}} \leq AL .
\]
This implies that $u$ has an endpoint $x_1\in\M$, and
\[
\lim_{t\to L} \frac1{\abs{h\circ u(t)}} = \Abs{(\pi\inner\epsilon)(x_1)} \geq B^{-1}.
\]
So $\pi(x_1)\neq 0$ and $x_1\in\Sigma$, thus $u$ can be continued through $t=L$, and therefore $(\Sigma,ds'^2)$ is geodesically complete.

Let $\tilde\Sigma$ be the universal covering of a connected component of $\Sigma$. This is geodesically complete, flat, and simply connected, therefore it is isometric to the Euclidean plane.
We can choose Cartesian coordinates $x$ and $y$ over $\tilde\Sigma$ such that,
\[
h = a + bx^2 + cy^2
\mbox.
\]
The Gaussian curvature of $ds^2$ is
\begin{align*}
\tfrac12 R &= h \,\nabla^2h - (\nabla h)^2\\
 &= 2a(b+c) - 2b (b-c) x^2 + 2c(b-c) y^2
\mbox.
\end{align*}
Of course, the curvature of a compact surface is bounded, so this must be bounded. In other words, the $x^2$ and $y^2$ terms must vanish, hence $c=b$.

There are now two possibilities. If $b=0$, then $ds^2$ is flat. Because $h=a$ is constant, it does not diverge even at infinity, hence $\pi$ does not vanish anywhere on $\M$, and $\Sigma=\M$. In this case $\M$ is a flat, compact Riemannian surface, therefore it is a torus, and $\pi = a^{-1}\epsilon^{-1}$ is constant.

If $b\neq0$, then it must have the same sign as $a$ (so that $h$ does not vanish). By the above formula, the curvature is the constant $4ab>0$, so $\M$ must be a sphere. The metric is
\[
ds^2 = (a + b[x^2+y^2])^{-2}(dx^2+dy^2)
\mbox.
\]
We can put this in a more standard form if we use the complex coordinate $\zeta := \left(\frac{a}{b}\right)^{1/2} (x + iy)$:
\[
ds^2 = (ab)^{-1}(1+\zeta\bar\zeta)^{-2} d\zeta\, d\bar\zeta
\mbox.
\]
The Poisson bracket is given by,
\[
\{\zeta,\bar\zeta\} = 2i \{x,y\} = 2i h = 2ai (1+ \zeta\bar\zeta)
\mbox.\]
This is the Poisson structure corresponding to the \Podles\ standard sphere.
\end{proof}

\section{Divergence}
\label{divergence.sec}
In \cite{haw} I derived $d(\pi\inner\epsilon)=0$ as a compatibility condition between a Poisson structure and volume form. The result of this section allows this condition to be restated in terms of a contravariant connection. I present this here partly because it facilitates the simplest proof of Lem.~\ref{volume.lift}, but I am stating a much more general theorem than I will really need, because it may be of independent interest.

Given a contravariant connection $\D$ and a differential form $\sigma\in\Omega^p(\M)$, we can define a ``contravariant divergence'' by $(\D\cdot\sigma)_{i_1i_2\dots i_{p-1}} := \D^i\alpha_{ji_1\dots i_{p-1}}$.
\begin{thm}
\label{D.divergence}
Let $\D$ be a torsion-free contravariant connection. There is a vector field $\phi\in\Gamma(\M,T\M)$ such that for any $\sigma\in\Omega^*(\M)$
\beq
\D\cdot\sigma = \phi\inner\sigma - \delta\sigma
\eeq
where $\delta\sigma := \pi\inner d\sigma - d(\pi\inner\sigma)$ is the Koszul-Brylinski codifferential \cite{bry,kos} used in Poisson homology. In particular, for $\alpha\in\Omega^1(\M)$
\beq
\label{D.divergence.eq}
\D\cdot \alpha = \phi\inner \alpha - \pi\inner d\alpha
\eeq
and for a volume form $\epsilon\in\Omega^n(\M)$,
\beq
\label{volume}
\phi\inner\epsilon = -d(\pi\inner\epsilon)
\eeq
if and only if $\D\epsilon=0$.
\end{thm}
\begin{proof}
We first note that for any $f,g\in\C^\infty(\M)$
\[
\D\cdot(f\,dg) = f\, \D\cdot dg + \{g,f\}
\mbox.\]
This gives the product identity 
\[
\D\cdot d(fg) = f\, \D\cdot dg + g\, \D\cdot df
\]
meaning that $\D\cdot d$ is a first order differential operator and is equivalent to some $\phi\in\Gamma(\M,T\M)$ as $\D\cdot df = \phi(f)$.

Now,
\[
\D\cdot (f\,dg) = \phi\inner (f\,dg) - \pi(df,dg)
\]
which implies eq.~\eqref{D.divergence.eq}.

Let $\sigma\in\Omega^p(\M)$ and $X\in\Gamma(\M,\Wedge^{p-1}T\M)$. The contravariant exterior derivative can be expressed as
\[
-[X,\pi] = \D\wedge X
\mbox.\]
 This is dual to the Koszul-Brylinski codifferential:
 \[
 [X,\pi]\inner\sigma = (-)^p X\inner\delta\sigma + \delta(X\inner\sigma)
 \mbox.\]
 We can also write this in terms of the contravariant divergence
 \begin{align*}
 (\D\wedge X)\inner\sigma 
 &= \D\cdot(X\inner\sigma) + (-)^p X\inner (\D\cdot\sigma) \\
 &= -(-)^pX\inner(\phi\inner\sigma) - \delta(X\inner\sigma) + (-)^p X\inner (\D\cdot\sigma)
 \mbox.
 \end{align*}
 From these we can solve for $X\inner(\D\cdot\sigma)$ and prove the general result.

The special case of a volume form $\epsilon\in\Omega^n(\M)$ follows from the simplification $\delta\epsilon = -d(\pi\inner\epsilon)$. Because $\epsilon\in\Omega^n(\M)$, $\D\cdot\epsilon=0$ if and only if $\D\epsilon=0$.
\end{proof}

The last result here shows that if $\D\epsilon=0$ then $\phi$ is a modular vector (see \cite{wei3}). The operator $\delta$ is the boundary operator defining Poisson homology. The modular class in Poisson homology is independent of $\epsilon$. It is easy to construct $\D$ to give $\phi$ any desired value, but unless $\phi$ belongs to the modular class, there cannot exist a volume form compatible with $\D$.

\section{Realizations}
\label{realizations}
Metacurvature is rather difficult to compute in general. It is only defined when $\D$ is flat and torsion-free, but these are highly nontrivial conditions. In principle, it is possible to write an explicit formula for the metacurvature of a metric contravariant connection, but this would be hopelessly complicated. Only in the symplectic case is the metacurvature easy to compute and understand. For this reason, my strategy for analyzing the condition $M=0$ on a Poisson manifold is to relate the Poisson manifold to a symplectic manifold. Fortunately, there is a well established way of doing this: symplectic realization.

\begin{definition}
A \emph{Poisson map} is a diffeomorphism of Poisson manifolds $\varphi : \M_1 \to \M_2$ such that the pull-back of functions intertwines the Poisson brackets,
\[
\varphi^*\{f,g\}_2 = \{\varphi^* f,\varphi^* g\}_1
\mbox.
\]
A (local) \emph{symplectic realization} of a Poisson manifold $\M$ is a Poisson map $\varphi : \bar\M \to \M$ with $\bar\M$ symplectic. A symplectic realization is \emph{full} if it is a submersion.
\end{definition}
There are variations on this definition in the literature \cite{wei1,vai2}. It is common to assume that $\varphi$ is a surjective submersion and a complete Poisson map. However, I will not need such properties, and I do not want to assume \emph{a priori} the global integrability necessary to satisfy them. Instead, I will use full local realizations as a tool for studying the local differential geometry of a Poisson manifold.

Given a symplectic realization, the preimages of points in $\M$ are the leaves of a foliation $\F$.  The symplectic orthogonal subbundle 
\[
T\Fp := \{v\in T\bar\M \mid \omega(v,w)=0 \; \forall w\in T\F\}
\] 
is integrable and thus defines another foliation, $\Fp$.

The symplectic foliation $\S$ of a Poisson manifold is defined by its tangent distribution. The tangent fiber over $x\in\M$ is
\[
T_x\S := \Im \#_x \subseteq T_x\M
\mbox.
\]
In general, the symplectic foliation is only a singular foliation, but over any open region where $\rk \pi$ is constant, $\S$ is a regular foliation and $T\S$ is actually a bundle. Define the set of regular points $\Mreg$ be the union of open sets over which $\rk \pi$ is constant.

\begin{definition}
An \emph{isotropic} realization is a full symplectic realization such that $\F\subset\Fp$.
\end{definition}
If $\varphi:\bar\M\to\M$ is an isotropic realization, then $\bar\M$ is of minimal dimension and $\Im\varphi\subseteq\Mreg$.
In this case, the symplectic leaves in $\M$ are the images of the $\Fp$-leaves, and  the $\Fp$-leaves are precisely the preimages of the symplectic leaves.

\subsection{Realizations and Connections}
Recall that for a foliated manifold, a foliated chart is a coordinate chart with two types of coordinates: transverse and leafwise. This is such that any leaf in the coordinate neighborhood is identified with a subset where the transverse coordinates are constant. The leafwise coordinates are thus coordinates along the leaf and the transverse coordinates are coordinates on the set of leaves.

On a symplectic manifold, a contravariant connection is equivalent to a covariant connection; they are intertwined by the map $\#$. In general, if $\M$ is a Poisson manifold with a symplectic realization $\varphi:\bar\M\to\M$, the composed map $\#\varphi^*:\Omega^1(\M)\to\Gamma(\bar\M,T\bar\M)$ plays the role that $\#$ does in the symplectic case. This gives a correspondence between the contravariant geometry of $\M$ and geometry on the leaves of $\Fp$.

Let $n=\dim\M$ and $2N=\dim\bar\M$. Thus $\dim\F=2N-n$ and $\dim\Fp=n$.

\begin{lem}
Let $\varphi:\bar\M\to \M$ be a full symplectic realization. If $\D$ is a contravariant connection on $\M$, then there exists a unique partial connection $\nabla$ on $\Fp$ such that, for $\alpha,\beta\in\Omega^1(\M)$,
\beq
\label{D.realization}
\#\varphi^*(\D_\alpha\beta) = \nabla_{\#\varphi^*\alpha} (\#\varphi^*\beta)
\mbox.\eeq
If $\D$ is flat or torsion-free, then so is $\nabla$. If $\D$ is flat and torsion-free then around any point of $\bar\M$, there exists an $\Fp$-foliated coordinate chart such that $\nabla$ is simply given by partial derivatives.
\end{lem}
\begin{proof}
First, observe that for $\alpha\in\Omega^1(\M)$, the pullback $\varphi^*\alpha\in\Omega^1(\bar\M)$ is normal to the distribution $T\F$. Therefore $\#\varphi^*\alpha\in \Gamma(\bar\M,T\Fp)$. Sections of this form span $\Gamma(\bar\M,T\Fp)$ as a $\C^\infty(\bar\M)$-module, therefore eq.~\eqref{D.realization} effectively defines $\nabla_X Y$ for all $X,Y\in\Gamma(\bar\M,T\Fp)$ because \eqref{D.realization} is consistent with the product rule.

The definition of a Poisson map, and the identity $d\{f,g\}=[df,dg]_\pi$ imply that $\varphi^*$ intertwines Koszul brackets. This then shows that for $\alpha,\beta\in\Omega^1(\M)$
\[
[\#\varphi^*\alpha,\#\varphi^*\beta] 
= \# [\varphi^*\alpha,\varphi^*\beta]_{\bar\pi}
= \#\varphi^*[\alpha,\beta]_\pi
\mbox.\]
This identity implies that the torsions and curvatures of $\D$ and $\nabla$ are intertwined by $\#\varphi^*$. Therefore if $\D$ is flat or torsion-free, then so is $\nabla$. 

Around any point of $\bar\M$, there exists a neighborhood $U\subseteq \bar\M$ such that the leaves of $\Fp|_U$ are simply connected and the leaf space $U/\Fp$ is Hausdorff and contractible. 
The flat, torsion-free partial connection $\nabla$ is precisely equivalent to a locally affine structure on each leaf of $\Fp$. Because each leaf of $\Fp|_U$ is simply connected, it can be identified as an open subset of an affine space (of dimension $n$). These form a bundle of affine spaces over $U/\Fp$. 
Because $U/\Fp$ is contractible, it can be identified with an open subset of $\R^{2N-n}$ and there exists a trivialization of the bundle of affine spaces. Together, this gives an identification of $U$ with an open subset of $\R^{2N-n}\times\R^n$. This is the desired foliated chart. 
\end{proof}

\begin{definition}
A \emph{cotangent curve} \cite{g-g} is a curve $(u,\xi): \R\supset I \to T^*\M$ such that $\#\xi(t) = \dot{u}(t)$. A \emph{cotangent geodesic} \cite{fer1} is a cotangent curve such that 
\beq
\label{geodesic.eq}
\D_\xi\xi=0
\mbox.
\eeq
\end{definition}

\begin{lem}
\label{geodesic}
Let $\varphi:\bar\M \to \M$ be a full symplectic realization, and $\D$ a contravariant connection on $\M$. 
\begin{enumerate}
\item
Any curve $v:I\to\bar\M$ in an $\Fp$-leaf descends to a unique cotangent curve $(u,\xi):I\to T^*\M$ such that $u=\varphi\circ v$ and $\dot v(t) = \#\varphi^*_{v(t)}(\xi(t))$. 
\item
For any cotangent curve $(u,\xi)$ and $t_0\in I$, there exists such an $\Fp$-curve $v$, defined over a neighborhood of $t_0$.
\item
 $(u,\xi)$ is a cotangent geodesic if and only if $v$ is a geodesic in an $\Fp$-leaf.
\end{enumerate}
\end{lem}
\begin{proof}
By assumption, $\varphi$ is a submersion. This implies that $\#\varphi^*_{v(t)}:T^*_{\varphi[v(t)]}\M \to T_{v(t)}\bar\M$ is injective. Its image is $\Im \#\varphi^*_{v(t)} = T_{v(t)}\Fp$.

To construct $(u,\xi)$ from $v$, let $u:=\varphi\circ v$; $\xi$ is uniquely defined by $\dot v(t) = \#\varphi^*_{v(t)}(\xi(t))$ because $\dot v(t)\in T\Fp$.

For the second claim, we must first choose a point in the preimage $\varphi^{-1}(u(t_0))$. Integrating $\dot v(t) = \#\varphi^*_{v(t)}(\xi(t))$ then defines $v(t)$ in a neighborhood of $t_0$.

The third claim follows from the definition of the lifted partial connection; the identity $\nabla_{\dot v}\dot v = \#\varphi^*(\D_\xi\xi)$ identifies the geodesic equations for $v$ and $(u,\xi)$.
\end{proof}

For the remainder of this section, $\D$ is a flat and torsion-free contravariant connection and $\varphi:\bar\M\to\M$ is a (local) full symplectic realization. 
\begin{definition}
A \emph{flat} $\Fp$-foliated chart is one where the induced partial connection $\nabla$ is given trivially by partial derivatives.
Let us say that a tensor on $\bar\M$ is $\Fp$-\emph{constant} (or \emph{linear}, or \emph{quadratic}, or \emph{polynomial}) if in any flat $\Fp$-foliated chart it is constant (linear, quadratic, polynomial) in the leafwise coordinates along $\Fp$.
\end{definition}

In fact, (although I won't need to prove it) a tensor is $\Fp$-polynomial if this is satisfied for \emph{some} flat $\Fp$-foliated atlas. These concepts really only depend upon the partial connection $\nabla$. The reason is that $\nabla$ extends to a flat partial connection on $T\bar\M$ with respect to $\Fp$, and this is unique modulo linear changes. This fact is related to the existence of a natural flat partial connection on the conormal bundle to a foliation. 

\begin{lem}
$M=0$ if and only if the Poisson bivector $\bar\pi$ of $\bar\M$ is $\Fp$-quadratic.
\end{lem}
\begin{proof}
Around any point of $\bar\M$, consider a flat foliated chart in a neighborhood $U$. With the symplectic structure of $\bar\M$, the trivial connection (partial derivatives) defines a flat, torsion-free contravariant connection on $U$. Consider the metacurvature $\bar M$ of $U$, and its relationship to the metacurvature $M$ of $\M$.

The definition of a Poisson map states that $\varphi^*$ intertwines Poisson brackets of functions. The relationship between $\D$ and $\nabla$ means that $\varphi^*:\Omega^*(\M)\to\Omega^*(U)$ intertwines the Poisson brackets of functions and $1$-forms. The product and Leibniz identities then imply that $\varphi^*$ intertwines all these generalized Poisson brackets of differential forms.

The definition of the metacurvature then shows that metacurvatures are intertwined as follows: For $\alpha,\beta,\gamma\in\Omega^1(\M)$,
\beq
\varphi^*\left[M(\alpha,\beta,\gamma)\right] = \bar M(\varphi^*\alpha,\varphi^*\beta,\varphi^*\gamma)
\mbox.\eeq

Now, assume that $M=0$. This implies that $0= \bar M(\alpha,\beta,\gamma)$ for any $\alpha,\beta,\gamma$ normal to $\F$. Because $U$ is symplectic, we can actually compute $\bar M$ from the third derivative of the Poisson bivector $\bar \pi$ on $U$. Equation \eqref{covariant.M} shows that the assumption $M=0$ is equivalent to the vanishing of all third derivatives of $\bar\pi$ with respect to leafwise coordinates in this chart.
In other words, $\bar\pi$ is $\Fp$-quadratic in this chart. However, since this works for any flat chart around any point, we can say that $\bar\pi$ is $\Fp$-quadratic.
\end{proof}

\begin{lem}
\label{volume.lift}
If there exists $\epsilon\in\Omega^n(\M)$ such that $d(\pi\inner\epsilon)=0$ and $\D\epsilon=0$, then the symplectic volume form $\frac1{N!}\omega^N$ is $\Fp$-constant.
\end{lem}
\begin{proof}
Again, around an arbitrary point of $\bar\M$, consider a flat $\Fp$-foliated chart. For any $\alpha\in\Omega^1(\M)$, the pull-back of its contravariant divergence is a divergence with the partial derivative connection $\partial$:
\[
\varphi^*(\D\cdot\alpha) = \partial\cdot(\#\varphi^*\alpha)
\mbox.
\]
 
In the notation of Thm.~\ref{D.divergence}, $d(\pi\inner\epsilon)=0$ and $\D\epsilon=0$ mean that $\phi=0$. By eq.~\eqref{D.divergence.eq}, for any $f\in\C^\infty(\M)$, $\D\cdot df =0$. Lifting this to $\bar\M$, we have 
\beq
0 = \varphi^*(\D\cdot df) 
= \partial\cdot (\#\varphi^* df) \label{partial.divergence}
\mbox.
\eeq
Now, 
\[
(\#\varphi^*df)\inner \tfrac{\omega^N}{N!} = d(\varphi^*f) \wedge \tfrac{\omega^{N-1}}{(N-1)!}
\]
is exact, so
\begin{align*}
0 &= d\left[(\#\varphi^*df)\inner \tfrac{\omega^N}{N!}\right] \\
&= \partial\cdot (\#\varphi^*df) \, \tfrac{\omega^N}{N!} + \partial_{\#\varphi^*df} \tfrac{\omega^N}{N!}
\mbox.
\end{align*}
The first term vanishes by eq.~\eqref{partial.divergence}. At any point of $\bar\M$, $\#\varphi^*df$ can give any vector in $T\Fp$. Therefore the symplectic volume form $\tfrac{\omega^N}{N!}$ is $\Fp$-constant.
\end{proof}

\begin{cor}
\label{symplectic.polynomial}
If $M=0$ and there exists $\epsilon\in\Omega^n(\M)$ with $\D\epsilon=0$ and $d(\pi\inner\epsilon)=0$, then the symplectic form on $\bar\M$ is $\Fp$-polynomial.
\end{cor}
\begin{proof}
\[
\omega = (-\bar\pi)^{N-1}\inner\frac{\omega^N}{N!}
\]
\end{proof}

\begin{lem}
\label{g.lift.thm}
If $\D$ is compatible with a metric $ds^2$ on $\M$, then there exists a unique flat metric $ds'^2$ (with the same signature) on the leaves of $\Fp$ such that,
\beq
\label{g.lift}
\langle\#\varphi^*\alpha,\#\varphi^*\beta\rangle' = \varphi^*[\langle\alpha,\beta\rangle]
\eeq
where $\langle\,\cdot\,,\,\cdot\,\rangle'$ is the $ds'^2$ inner product on $T\Fp$ and  $\langle\,\cdot\,,\,\cdot\,\rangle$ is the $ds^2$ inner product on $T^*\M$.
\end{lem}
\begin{proof}
This expression clearly defines an inner product on every fiber of $T\Fp$. As above, $\D$ determines a flat partial connection $\nabla$ along $\Fp$. Because of the way the metrics and connections are intertwined, the $\nabla$ derivative of the lifted metric $ds'^2$ vanishes. 

Along any leaf of $\Fp$, $\nabla$ becomes simply a (flat) connection, and $ds'^2$ a metric. This connection is the Levi-Civita connection, thus the metric is flat.
\end{proof}

\section{Riemannian Manifolds}
\label{Riemannian}
A noncommutative deformation of differential forms is characterized by a contravariant connection for which the torsion, curvature, and metacurvature vanish. In \cite{haw} I showed that if a Riemannian manifold is deformed into a real spectral triple, then the contravariant connection is compatible with the metric. I will now briefly present a more robust derivation of this condition, using a much weaker notion of spectral triple.  I expect that this compatibility condition can be derived from any concrete notion of noncommutative geometry.

Recall that a spectral triple $(\A,\Hi,D)$ consists of a Hilbert space $\Hi$, an involutive algebra of bounded operators $\A$, and a self-adjoint unbounded operator $D$, such that the commutator of $D$ with any element of $\A$ is bounded. 

The geometry of a Riemannian manifold $\M$ can be encoded algebraically in a spectral triple. Let $\A=\Ci_0(\M)$, let $D$ be any Dirac-type operator, and let $\Hi$ be the Hilbert space of square-integrable sections of the bundle on which $D$ acts. The metric can be recovered because of the identity
\[
([D,f])^2 = - \langle df,df\rangle
\]
where $f\in\Ci_0(\M)$ and $\langle\,\cdot\,,\,\cdot\,\rangle$ is the metric pairing. A spectral triple is in this sense a generalization of a Riemannian manifold.

Differential forms can be constructed from a spectral triple (see \cite{con1} and Sec.~\ref{forms}). In particular, $\Omega^1_D(\A)$ is the $\A$-bimodule generated by bounded operators of the form $da:=[D,a]$, for $a\in\A$. In the above example, $\Omega^1_D(\A)\cong\Omega^1_0(\M)$.

Now suppose that there exists a smooth noncommutative deformation of a spectral triple describing a Riemannian manifold. Suppose that the differential forms constructed from this deformation are a smooth deformation of the differential graded algebra $\Omega^*_0(\M)$. We have seen in Section \ref{meta} that this deformation of $\Omega^*_0(\M)$ is described to leading order by a Poisson structure and a contravariant connection $\D$ with $0$ torsion, curvature, and metacurvature.

To see how this contravariant connection is related to the metric, consider the simple identity 
\beq
\label{simple.identity}
[a,([D,b])^2] = [a,[D,b]]\,[D,b] +  [D,b]\,[a,[D,b]]
\mbox.
\eeq
A commutator such as $[a,[D,b]]$ corresponds to a generalized Poisson bracket such as $\{f,dh\}$. At first order in $\k$, the identity \eqref{simple.identity} gives
\[
\{f,\langle dh,dh\rangle\} = 2 \langle dh,\{f,dh\}\rangle = 2 \langle dh,\D_{df}dh\rangle
\mbox,
\]
but this is simply the condition that the contravariant derivative of the metric is $0$. Since $\D$ is torsion-free, this means that $\D$ is (by definition) the metric contravariant connection.

Together with the compatibility with the volume form, this gives the following notion of compatibility.
\begin{definition}
\label{compatible}
For the purposes of this paper, I will say that a metric and Poisson structure are \emph{compatible} if
\begin{enumerate}
\item
The metric contravariant connection $\D$ is flat.
\item 
The metacurvature (of $\D$) vanishes: $M=0$.
\item
the Poisson structure is compatible with the Riemannian volume form: 
\[
d(\pi\inner\epsilon)=0
\mbox.
\]
\end{enumerate}
\end{definition}

\subsection{Local Structure}
\label{local}
In this subsection, $\M$ is a compact, connected Riemannian manifold with a metric tensor $g_{ij}$ and a compatible Poisson structure $\pi$. The analysis here is based on repeatedly applying two simple principles: A continuous function on a compact manifold is bounded, and a bounded polynomial is constant.

The local norm of the Poisson bivector $\pi$ is defined by
\[
\abs{\pi}^2 := \tfrac12 \pi^{ij}\pi^{kl}g_{ik}g_{jl}
\mbox.\]
This is a continuous function and $\M$ is compact, so there is also a global norm,
\[
\norm{\pi} := \max_{x\in\M} \abs{\pi(x)}
\mbox.
\]
The norm of any multivector or differential form is defined analogously.

Dual to the inclusion $T\Fp\subset T\bar\M$, there is a restriction map $\Omega^*(\M)\to\Omega^*(\Fp)$. Locally, this gives the pull-back of a differential form to every leaf of $\Fp$. Let $\omega_0\in\Omega^2(\Fp)$ be the restriction of the symplectic form $\omega\in\Omega^2(\bar\M)$.

\begin{lem}
For any local full realization of $\M$,  $\omega_0$ is $\Fp$-constant.
\end{lem}
\begin{proof}
The cotangent geodesic flow is given by a smooth vector field on $T^*\M$. The geodesic equation \eqref{geodesic.eq} preserves the norm $\abs{\xi(t)}$, so the geodesic flow is tangent to the sphere subbundles 
\[
\left\{\xi\in T^*\M \bigm| \abs\xi = \mathrm{const.}\right\}
\]
which are compact. This is thus a complete vector field and integrates to give a complete cotangent geodesic $(u,\xi):\R\to T^*\M$ through any point in $T^*\M$.

For any $t_0\in \R$, there exists a symplectic realization $\varphi: \bar\M\to \M$ over $u(t_0)$. By Lem.~\ref{g.lift.thm}, this defines flat Riemannian metrics on the $\Fp$-leaves in $\bar\M$. By Lem.~\ref{geodesic}, $(u,\xi)$ lifts to a geodesic $v$ in the $\Fp$-leaf over $u$. The pull-back of the local norm $\abs{\pi}$ to $\bar\M$ is 
\[
\varphi^*\abs{\pi} = \abs{\omega_0}
\mbox,
\]
the norm of $\omega_0$ with respect to the flat metric on $\Fp$-leaves. By Cor.~\ref{symplectic.polynomial}, $\omega$, and thus $\omega_0$ and $\abs{\omega_0}^2$ are $\Fp$-polynomial. This means that 
\[
\abs{\pi(u(t))}^2 = \abs{\omega_0(v(t))}^2
\]
is a polynomial for $t\approx t_0$. Hence, $\abs{\pi(u(t))}^2$ is a polynomial for all $t\in \R$.

This polynomial is bounded by $\norm{\pi}^2$, therefore it is constant. Through every point of $\M$ there are cotangent geodesics passing in any direction tangent  to the symplectic foliation $\S$, therefore $\abs{\pi}$ is $\S$-constant.

Now, for any symplectic realization, this shows that $\abs{\omega_0}=\varphi^*\abs{\pi}$ is $\Fp$-constant. Because $\omega_0$ is $\Fp$-polynomial, this is enough to prove that $\omega_0$ is $\Fp$-constant.
\end{proof}

\begin{lem}
Around any regular point $x_0\in\Mreg$ (with $\rk\pi(x_0)=:2m$) there exists an $\S$-foliated chart with leafwise coordinates $\{x^a\}_{a=1}^{2m}$ and transverse coordinates $\{z^\alpha\}_{\alpha=1}^{n-2m}$ such that:
\begin{enumerate}
\item
The components of $\pi$ are constant in this coordinate system.
\item
The metric is of the form
\beq
\label{reduced.g}
ds^2 = g^{\scriptscriptstyle\parallel}_{ab}(dx^a - A^a_\alpha dz^\alpha)(dx^b - A^b_\beta dz^\beta) +  g^\perp_{\alpha\beta} dz^\alpha dz^\beta
\eeq
with $g^{\scriptscriptstyle\parallel}_{ab}$ and $g^\perp_{\alpha\beta}$ independent of $x$. $A^a_\alpha$ is polynomial in $x$, and for fixed $z$ and $\alpha$, $A^a_\alpha$ are the components of a Hamiltonian vector field in the symplectic leaf.
\end{enumerate}
\end{lem}
\begin{proof}
First consider any $\S$-constant function $Z$, defined over a neighborhood of $x_0$. This means that $\# dZ=0$. A direct calculation shows that $0 = \D_i dZ_j + \D_j dZ_i$ (see \cite[eq.~(6.4)]{haw}). This is the analogue of Killing's equation. If $\varphi$ is an isotropic realization, then $\#d(\varphi^*Z)$ is tangent to $\Fp$ and restricts to a Killing vector in each $\Fp$ leaf. Any two such functions Poisson-commute and thus define commuting Killing vectors. Note that $\varphi^*Z$ is $\Fp$-constant.

Let $\varphi:\bar\M\to\M$ be an isotropic realization over $x_0$, with $\bar\M$ small enough that the leaf space $\bar\M/\Fp$ is Hausdorff, connected, and simply connected.   
Let $\{z^\alpha\}_{\alpha=1}^{n-2m}$ be a system of coordinates on $\bar\M/\Fp$, or the equivalent $\Fp$-constant functions on $\bar\M$. 

Define $\varepsilon_{(\alpha)}:= \#dz^\alpha\in\Gamma(\bar\M,T\F)$ for $\alpha=1,\dots,n-2m$. These are mutually commuting Killing vectors in the leaves of  $\Fp$. They form a basis for $T\F$ because 
\[
\varepsilon_{(1)}\wedge\dots\wedge\varepsilon_{(n-2m)} = \#\left(dz^1\wedge\dots\wedge dz^{n-2m}\right)
\]
is nonvanishing. 

Let $\mathcal T\subset \bar\M$ be a smooth transversal to $\Fp$ (intersecting $\varphi^{-1}(x_0)$); this is parametrized by the $z$'s. Choose $\{e_{(a)}\}_{a=1}^{\;2m} \subset \Gamma(\mathcal T, T\Fp)$, a basis of the orthogonal complement of $T\F\subset T\Fp$ (along $\mathcal T$) such that $\omega(e_{(a)},e_{(b)})$ are constant. Extend these vectors to a neighborhood of $\mathcal T$ by requiring that they be $\Fp$-constant: $0=\nabla e_{(a)}$. These are (in particular) mutually commuting Killing vectors in the leaves of $\Fp$.

Because $\omega_0$ is $\Fp$-constant, the partial connection $\nabla$ preserves $T\F = \ker \omega_0 \subset T\Fp$ and the vectors $\varepsilon_{(\alpha)}$ are orthogonal to $T\F$ everywhere. Because $\{\varepsilon_{(\alpha)}\}_{\alpha =1}^{\,n-2m}$ are Killing vectors spanning $T\F$, any derivative orthogonal to $T\F$ must vanish, so
\[
[e_{(a)},\varepsilon_{(\alpha)}] = \nabla_{e_{(a)}}\varepsilon_{(\alpha)} - \nabla_{\varepsilon_{(\alpha)}}e_{(a)} = 0 .
\]

So, we have a basis of mutually commuting vectors tangent to $\Fp$ near $\mathcal T$. Exponentiating these, we can construct coordinates $\{x^a\}_{a=1}^{2m}$, $\{y^\alpha\}_{\alpha=1}^{n-2m}$, and $\{z^\alpha\}_{\alpha=1}^{n-2m}$ in a neighborhood of $\mathcal T$ such that $\frac{\partial}{\partial x^a} = e_{(a)}$ and $\frac{\partial}{\partial y^{\alpha}} = \varepsilon_{(\alpha)}$.

This is an $\Fp$-foliated chart with transverse coordinates $z^\alpha$. It is also an $\F$-foliated chart with transverse coordinates $x^a$ and $z^\alpha$. To put this another way, the $z$'s parametrize the set of $\Fp$-leaves, the $x$'s parametrize the set of $\F$-leaves in each $\Fp$-leaf, and the $y$'s are coordinates on each $\F$-leaf, as depicted in the figure.
\begin{figure}
\includegraphics{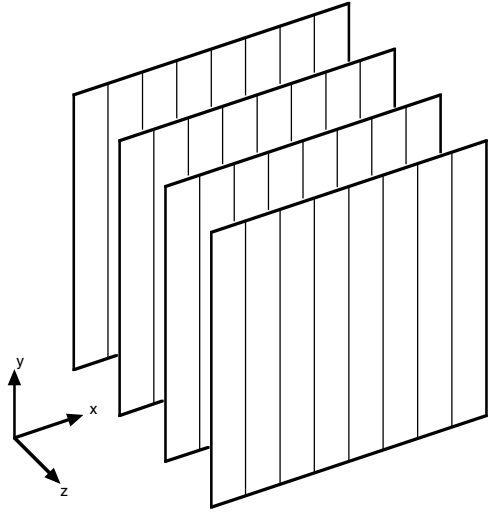}
\caption{In this sketch, the planes represent $\Fp$-leaves; the lines represent $\F$-leaves.}
\end{figure}

Because the $\varepsilon_{(\alpha)}$'s are orthogonal to the $e_{(a)}$'s, the $\Fp$ leaf metric takes the form,
\[
ds'^2 = g'_{ab}dx^a dx^b + g'_{\alpha\beta} dy^\alpha dy^\beta
\mbox.\]
Because these are Killing vectors, the components $g'_{ab}$ and $g'_{\alpha\beta}$ are functions of the $z$'s alone. This shows that this is a \emph{flat} $\Fp$-foliated chart.

The $y$-$y$-components of $\omega$ vanish because $\F$ is isotropic. The $x$-$y$-components vanish because (by definition) $\Fp$ is symplectically orthogonal to $\F$. By definition, $\varepsilon_{(\alpha)}\inner\omega = dz^{\alpha}$.
With this in mind, the symplectic form reduces to
\[
\omega = \tfrac12 \omega_{ab} dx^a\wedge dx^b + \omega_{a\beta}dx^a\wedge dz^\beta+ dy^\alpha\wedge dz^\alpha
\mbox.\]
The first term is $\omega_0$, and therefore $\omega_{ab}=\omega(e_{(a)},e_{(b)})$ is constant.

The Poisson bivector $\bar\pi$ on $\bar\M$ is of course given by inverting this symplectic form. The Poisson bivector on $\Im\varphi\subset\M$ is given by the $x$ and $y$ components of $\bar\pi$, but of those, only the $x$-$x$-components may be nonzero. Thus the nonzero part of $\pi$ is $\pi^{ab}$ which is the inverse matrix of $\omega_{ab}$ and is constant in these coordinates.

The metric on $\Im\varphi\subset\M$ is constructed by using $\omega$ to lower the indices of the contravariant form of $ds'^2$. This gives
\[
ds^2 = g'^{cd}(\omega_{ca}dx^a + \omega_{c\alpha}dz^\alpha) (\omega_{db}dx^b+\omega_{d\beta}dz^\beta) + g'^{\alpha\beta} dz^\alpha dz^\beta
\mbox.\]
This is of the form \eqref{reduced.g} with
\[
A^a_\alpha :=  \pi^{ba}\omega_{b\alpha}
\mbox.\]
Because $\omega$ is closed and $\omega_{ab}$ constant, $\omega_{a\alpha}dx^a$ is exact. Therefore $A^a_\alpha$ is a Hamiltonian vector field in the sense claimed.
\end{proof}

These results can be restated in purely geometric terms:
\begin{enumerate}
\item
The symplectic foliation $\S$ is a Riemannian foliation \cite{mol} of $\Mreg$. That is, the metric descends to a well-defined transverse metric $g^\perp_{\alpha\beta}$. This is a metric on the leaf space, to the extent that the leaf space is meaningful.
\item
The induced leaf metric $g^{\scriptscriptstyle\parallel}_{ab}$ is flat, thus the leaves are locally affine and in this sense the symplectic form is constant on each leaf.
\item 
Locally, $\S$ looks like a bundle of symplectic affine spaces. Flowing orthogonally to $\S$ defines a connection $\lcon$ (for which $A^a_\alpha$ is the potential). The structure Lie algebra is that of polynomial Hamiltonian vector fields on $\R^{2m}$.
\end{enumerate}

\begin{cor}
The regular symplectic leaves are (intrinsically) geodesically complete.
\end{cor}
\begin{proof}
We can compute $\abs{\pi^m}$ explicitly from the above coordinate expressions. It is $\S$-constant. 

Suppose that $u:[0,1)\to\M$ is an incomplete, inextensible geodesic in a regular $2m$-di\-men\-sional  leaf. Because $\M$ is compact and $\abs{\dot u}$ is constant, this extends continuously to $u:[0,1]\to\M$. The point $u(1)$ does not lie on the leaf, thus $\rk \pi[u(1)] < 2m$ and $\abs{\pi^m[u(1)]}=0$. However, $\abs{\pi^m[u(t)]}$ is constant and nonzero for $t\in[0,1)$, contradicting the continuity of $\abs{\pi^m}$. Therefore there does not exist an incomplete, inextensible geodesic.
\end{proof}

The extrinsic curvature of $\S$ is $K := \frac12\lcon ds_{\scriptscriptstyle\parallel}^2$. This is a section of $T^*\S\otimes T^*\S\otimes N^*\S$, where $N^*\S$ is the conormal bundle to $\S$. This is easily computed for the metric form \eqref{reduced.g}:
\[
K_{ab\alpha} = \tfrac12\left(g^{\scriptscriptstyle\parallel}_{ab,\alpha} + A_{a\alpha,b} + A_{b\alpha,a}\right)
\]
where Latin indices are lowered with the leaf metric $g^{\scriptscriptstyle\parallel}_{ab}$. Commas denote partial derivatives.

\begin{lem}
$A^a_\alpha$ depends only linearly on the $x$'s (i.e., it is $\S$-linear).
\end{lem}
\begin{proof}
First, consider the trace of the extrinsic curvature:
\[
(\tr K)_\alpha := 
K^a_{\; a\alpha} = \tfrac12 g_{\scriptscriptstyle\parallel}^{ab} g^{\scriptscriptstyle\parallel}_{ab,\alpha} + A^a_{\alpha,a}
\mbox.\]
The last term vanishes because $A^a_\alpha$ is Hamiltonian. So $K^a_{\; a\alpha}$ is constant along $\S$. (This is actually a coordinate-independent statement; the foliation determines a canonical flat partial connection of $N^*\S$ along $\S$.)

The  components $R_{abcd}$ of the Riemann tensor parallel to $\S$ can be computed from the extrinsic curvature:
\[
R_{abcd} = K_{ac\alpha}K_{bd}^{\;\;\alpha} - K_{ad\alpha}K_{bc}^{\;\;\alpha}
\mbox.\]
The trace of this is
\beq
\label{R.trace}
g_{\scriptscriptstyle\parallel}^{ac}g_{\scriptscriptstyle\parallel}^{bd} R_{abcd} =
\tr K^2 - (\tr K)^2 = K_{ab\alpha}K^{ab\alpha} - K^a_{\;a\alpha}K^{b\;\alpha}_{\;b}
\mbox.
\eeq
The last term is $\S$-constant. This does not necessarily extend to a continuous function beyond $\Mreg$. However, \eqref{R.trace} is bounded in terms of the norm of the Riemannian curvature of $\M$. This means that $K$ must be bounded. Since $A$ is polynomial in $x$, so is $K$ and it must be independent of $x$.

This shows that $A_{a\alpha,b}+A_{b\alpha,a}$ is independent of $x$, which implies that $A^a_\alpha$ is linear in $x$.
\end{proof}

This shows that the structure group of $\lcon$ reduces to $\mathrm{ISp}(2m,\R)$, the group\footnote{``Inhomogeneous symplectic group''} of affine symplectomorphisms of $\R^{2m}$. This is very much like the geometry of a Kaluza-Klein model (see the review \cite{wit}). In a Kaluza-Klein model, Yang-Mills (gauge) theory is realized geometrically by a bundle of homogeneous spaces of the gauge group. The Yang-Mills Lagrangian (gauge curvature squared) is recovered as a term of the Einstein-Hilbert Lagrangian (Riemannian scalar curvature). Here, $\mathrm{ISp}(2m,\R)$ plays the role of the gauge group and $\lcon$ the gauge connection.

The curvature of $\lcon$ is
\[
F^a_{\alpha\beta} = A^a_{\beta,\alpha} - A^a_{\alpha,\beta} + A^b_\alpha A^a_{\beta,b} - A^b_\beta A^a_{\alpha,b}
\mbox.\]
This should be thought of as a vector field (over $\R^{2m}$) valued $2$-form (over the leaf space).

\begin{lem}
\label{reg.pi.form}
In some neighborhood of any regular point where $\rk\pi=2m$, there exist $2m$ commuting Killing vectors $\{X_A\}_{A=1}^{\;2m}$ which span $T\S$. The Poisson bivector can be expressed as
\beq
\label{pi.form}
\pi = \tfrac12 \Pi^{AB} X_A\wedge X_B
\eeq
where the matrix $\Pi^{AB}$ is constant and nondegenerate.
\end{lem}
\begin{proof}
The scalar curvature $R$ of $\M$ can be computed in terms of the transverse metric, the leaf metric, and $\lcon$. The only term of $R$ which is not necessarily $\S$-constant is the Yang-Mills term
\[
-\tfrac14 F^2 := 
-\tfrac14 g^{\scriptscriptstyle\parallel}_{ab} g_{\scriptscriptstyle\perp}^{\alpha\gamma} g_{\scriptscriptstyle\perp}^{\beta\delta} F^a_{\alpha\beta} F^b_{\gamma\delta}
\mbox.
\]
 Because $R$ is bounded, $F^a_{\alpha\beta}$ must be bounded. However, $F^a_{\alpha\beta}$ must be polynomial in $x$, therefore it is independent of $x$.

In other words, the curvature is translation valued. This means that locally, by a gauge (coordinate) transformation, $A^a_\alpha$ can be made independent of $x$. With such a coordinate choice, all components of the metric \eqref{reduced.g} are independent of $x$. This means that the basis vectors in $x$-directions are commuting Killing vectors. Let $\{X_A\}_{A=1}^{\;2m}$ be this basis; that is,
$X^a_A=\delta^a_A$ and $X^\alpha_A=0$. The components of $\pi$ in this coordinate system become the components of $\pi$ in this basis. Renaming these as $\Pi^{AB}$, we have eq.~\eqref{pi.form}.
\end{proof}

The final step is to show that Lem.~\ref{reg.pi.form} implies that the decomposition \eqref{pi.form} exists in a neighborhood of any point, not just a regular point. To do this, I will view Killing vectors as coming from a larger bundle. Let $\nabla$ be the Levi-Civita connection on $\M$.  If $X\in\G(\M,T\M)$ is a Killing vector, then $\nabla X$ is antisymmetric (with respect to the metric), so $X$ and $\nabla X$ together form a section of the bundle $\iso(T\M)\cong T\M \oplus \Wedge^2T\M$; moreover, the second derivative can be expressed in terms of $X$ and the Riemann tensor. 

Based on this, define a connection on $\iso(T\M)$ by
\[
\tnabla_Y(X,w) := \left(\nabla_YX-w(Y),\nabla_Yw - R(Y,X)\right)
\]
where $Y\in\G(\M,T\M)$ and $(X,w)\in\G(\M,\iso(T\M))$. If $X$ is a Killing vector then $0=\tnabla(X,\nabla X)$. Conversely, any $\tnabla$-constant section is given by a Killing vector in this way.

Let $2m$ be the maximum rank of $\pi$.
\begin{thm}
\label{local.structure}
Over any simply connected open subset of $\M$, $\pi$ is given by $2m$ commuting Killing vectors as 
\[
\pi = \tfrac12 \Pi^{AB} X_A\wedge X_B
\]
with constant coefficients. The set $\{(X_A,\nabla X_A)\}_{A=1}^{\;2m}$ is a basis of $\tnabla$-constant sections of a flat subbundle $V\subset\iso(T\M)$.
\end{thm}
\begin{proof}
This proof consists of repeatedly applying another simple principle: If a continuous section vanishes over some neighborhood of any regular point, then it vanishes over $\Mreg$ and by continuity it vanishes over $\M$. This can be applied to any property that can be expressed as the vanishing of a continuous section.

Let $x_0\in\Mreg$ be an arbitrary regular point. By Lem.~\ref{reg.pi.form}, in some neighborhood $U$ of $x_0$, we have the decomposition \eqref{pi.form}. From this, we can construct $\tilde X_A :=(X_A,\nabla X_A)$ and 
\beq
\label{pit.form}
\tilde\pi = \tfrac12\Pi^{AB} \tilde X_A\wedge \tilde X_B
\eeq
over $U$. Because each $X_A$ is a Killing vector, $0=\tnabla\tilde X_A$, and so
\[
0 = \tnabla\tilde\pi
\mbox.\]
This $\tilde\pi\in\G(U,\Wedge^2\iso(T\M))$ is made up of three components; the first is $\pi$ itself, and the others are sections of $T\M\otimes\so(T\M)$ and $\Wedge^2\so(T\M)$. 

\emph{A priori}, $\tilde\pi$ is only defined over $U$, and appears to depend upon a choice of decomposition. However, observe that
\[
A^{ij}_{\,k} := \tfrac12 \left(\pi^{ij}_{\;|k} - \pi^{j\;\: |i}_{\;k} - \pi^{i\;\: |j}_{\;k}\right)
= \Pi^{AB}X_A^i X^j_{B|k} 
\]
and
\[
B^{ij}_{kl} := A^{ij}_{\,l|k} + \pi^{im}R^j_{\: klm}
= \Pi^{AB} X^i_{A|k} X^j_{B|l}
\]
are the other components of $\tilde\pi$. Obviously, $A$ and $B$ are well defined tensors over $\M$. Although $B$ is not explicitly a section of $\Wedge\so(T\M)$, it is over $U$; by continuity, it is over $\M$.

So, we can define $\tilde\pi\in\G(\M,\Wedge^2\iso(T\M))$ as the global section with components $\pi$, $A$, and $B$. This satisfies $0=\tnabla\pi$ over $U$; by continuity, this is true over $\M$.

This implies that $\tilde\pi$ has constant rank (which must be $2m$). It thus spans a subbundle $V\subset \iso(T\M)$. Over $U$, the decomposition \eqref{pit.form} shows that the sections $\{\tilde X_A\}_{A=1}^{\;2m}$ span $V$. Since $0=\tnabla\tilde X_A$, the restriction of $\tnabla$ to $V$ is flat over $U$; by continuity, it is flat over $\M$.

Any $\tnabla$-constant section of $V$ over $U$ is a constant linear combination of $\{\tilde X_A\}_{A=1}^{\;2m}$; its $T\M$-component is the same linear combination of $\{X_A\}_{A=1}^{\;2m}$. Therefore, the $T\M$-components of $\tnabla$-constant sections of $V$ over $U$ are mutually commuting Killing vectors; by continuity, this is true over any domain in $\M$.

Now forget the $X_A$ and $\tilde X_A$ used above.
Over any simply connected neighborhood of any point in $\M$, there exists a basis $\{\tilde X_A\}_{A=1}^{\;2m}$ of $\tnabla$-constant sections of $V$. In this basis, $\tilde\pi$ has constant components, $\Pi^{AB}$. Defining $X_A$ as the $T\M$-component of $\tilde X_A$, this gives the desired decomposition of $\pi$.
\end{proof}

\subsection{Global Structure}
\label{global}
\begin{thm}
\label{structure}
Let $\M$ be a connected, compact Riemannian manifold with a compatible Poisson structure.
There exists a covering $\Mt$ and a Lie group $G$ of isometries of $\Mt$ such that: 
\begin{enumerate}
\item
$\M$ is the quotient $\Mt/\Gamma$ by a discrete, cocompact subgroup $\Gamma \subset G$.
\item
The Poisson structure on $\Mt$ is induced by an $\Ad_G$-invariant bivector $\Pi\in\Wedge^2\g$, where $\g$ is the Lie algebra of $G$.
\item
The span of $\Pi$ (in $\g$) densely generates a connected abelian normal subgroup $T \normal G$.
\item
$\G\cap T = \{e\}$ and the subgroup $T\,\G\subset G$ generated by $\G$ and $T$ is dense.
\end{enumerate}
\end{thm}
\begin{proof}
Choose some (arbitrary) base point in $\M$. Let $\G$ be the holonomy group for the flat bundle $V\subset\iso(T\M)$, regarded as a discrete group. $\G$ is a quotient of the fundamental group, so we can define $\Mt$ as the covering of $\M$ with covering group $\G$. Let $G_1$ be the (Lie) group of isometries of $\Mt$ that preserve the Poisson structure (i.e., Poisson isometries).

By construction, $\G\subset G_1$, $\M\cong\Mt/\G$, and $V$ is globally flat. So, the Killing vectors $X_A$ exist globally over $\Mt$, and $\pi=\tfrac12\Pi^{AB}X_A\wedge X_B$ over $\Mt$.

Since the Killing vectors $X_A$ commute, the decomposition shows that they preserve $\pi$. Hence, they are elements of the Lie algebra of $G_1$. 
Let $T\subset G_1$ be the Lie subgroup densely generated by $\{X_A\}_{A=1}^{\;2m}$. This is abelian because the $X_A$'s commute. Because $\pi$ is $G_1$-invariant, $\Pi$ is $\Ad_{G_1}$-invariant and $T$ is normal.

Define $G$ as the closure of the subgroup $T\,\Gamma \subset G_1$. Because $G\subseteq G_1$ is a closed subgroup, it is a Lie group. By construction $\Gamma\subset G$ and $T \normal G$.

The definition of $\G$ implies that the adjoint action of $\G$ on $T$ is a faithful representation. Because $T$ is abelian, the vectors $X_A$ are $T$-invariant. Therefore $\G\cap T = \{e\}$.

For some arbitrary point $x\in\Mt$, let $H_x \subset G \subset \Isom(\Mt)$ be the subgroup leaving $x$ fixed. Because $\M$ (and hence $\Mt$) is Riemannian (with positive definite metric) $H_x$ is compact. Let $\Or_x$ be the closure of the symplectic leaf through the image of $x$ in $\M$. This is naturally identified with the double quotient $H_x\backslash G /\Gamma$. Because $\M$ is compact, $\Or_x$ must be compact and so $G/\Gamma$ is compact. In other words, $\Gamma \subset G$ is cocompact.
\end{proof}
Note that an $\Ad_G$-invariant bivector in $\Wedge^2\g$ is the same thing as a bi-invariant (left and right invariant) bivector field on $G$. In fact \cite[Thm.~10.4]{vai2} a bi-invariant Poisson structure always comes from an abelian normal subgroup in this way.

\begin{thm}
Let $\M$ be an (arbitrary) Riemannian or pseudo-Riemannian manifold. If $\Mt$ is some covering of $\M$, and $G$ is a Lie group acting by isometries on $\Mt$ such that $\M$ is the quotient of $\Mt$ by a subgroup of $G$, then any $\Ad_G$-invariant Poisson bivector $\Pi\in \Wedge^2\g$ induces a Poisson structure on $\M$ which is compatible with the metric in the sense of Definition \ref{compatible}.
\end{thm}
\begin{proof}
Obviously, $\Pi$ induces a $G$-invariant Poisson structure on $\Mt$, but $\M$ is the quotient of $\Mt$ by a subgroup of $G$, so this induces a Poisson structure on $\M$. The compatibility conditions in question are all local, so it is sufficient to check them on $\Mt$.

By Theorem 10.4 of \cite{vai2}, $\Pi$ spans an abelian ideal $\ta\subseteq \g$. So $\Pi\in\Wedge^2\ta$. This means that the Poisson structure on $\Mt$ can be written as $\frac12\Pi^{AB}X_A\wedge X_B$ where $\{X_A\}$ is a basis of commuting Killing vectors spanning $\ta$.

The volume condition is straightforward to check ($\Lie$ denotes the Lie derivative):
\[
\begin{split}
d(\pi\inner\epsilon) &= \tfrac12\Pi^{AB} d[(X_A\wedge X_B)\inner\epsilon) \\
&= \Pi^{AB} X_A\inner (\Lie_{X_B}\epsilon) - \tfrac12 \Pi^{AB} [X_A,X_B]\inner\epsilon =0
\end{split}
\]
because the vectors $X_A$ preserve the volume form and commute.

Some contravariant connection $\D$ is defined by 
\[
\D_{\alpha}\beta = \Pi^{AB} (X_A\inner\alpha) \Lie_{X_B}\beta
\mbox,\]
for any $\alpha,\beta\in\Omega^1(\Mt)$.
This is compatible with the metric because the $X_A$'s are Killing vectors. 

It is sufficient to compute the torsion using exact $1$-forms. Noting that,
\[
\D_{df} dh = \Pi^{AB} X_A(f) d[X_B(h)]
\mbox,
\]
we have
\[
\begin{split}
T(df,dh) &= \D_{df}dh - \D_{dh}df - d\{f,h\} \\
&= \Pi^{AB}[X_A(f) d(X_Bh) + d(X_Af)X_B(h) - d(X_Af\, X_Bh)] \\
&=0\mbox.
\end{split}
\]
This shows that this $\D$ is in fact the metric contravariant connection.

Again, let $\ta$ be the abelian Lie algebra spanned by $\{X_A\}$. Given any regular point $x_0\in\Mt$, any covector at $x_0$ extends to a $\ta$-invariant $1$-form in a neighborhood of $x_0$. Let $\gamma$ be such a $\ta$-invariant $1$-form. Its contravariant derivative vanishes, $\D\gamma=0$. Consequently, $K(\alpha,\beta)\gamma=0$ and so $K=0$ at $x_0$. The curvature vanishes at every regular point, therefore $K=0$ everywhere.

The first expression for $\D$ generalizes to give the Poisson bracket of a function and a differential form,
\[
\{f,\sigma\} = \Pi^{AB} X_A(f) \Lie_{X_B}\sigma
\mbox.\]
The Leibniz identity \eqref{Leibniz2} then implies that the bracket of two $1$-forms must be,
\[
\{\alpha,\beta\} = \Pi^{AB} \Lie_{X_A}\alpha \wedge \Lie_{X_B}\beta
\mbox.
\]

Again, let $\gamma$ be a $\ta$-invariant $1$-form on some neighborhood of $x_0$. These formul\ae\ show that any generalized Poisson bracket of $\gamma$ with a function or $1$-form must vanish. Thus by eq.~\eqref{Mdef}, $M(df,\beta,\gamma) =0$. This shows that $M=0$ at $x_0$, hence at every regular point, hence everywhere.
\end{proof}
These two results show that in the case of a compact Riemannian manifold, a Poisson structure is compatible with the metric if and only if it is induced in this way from a bi-invariant Poisson structure on a group.

Theorem \ref{structure} thus shows how to construct all examples of compatible Poisson structures. They are classified by triples $(G,\Gamma,\Pi)$ of a Lie group, a cocompact discrete subgroup and a bi-invariant Poisson structure, such that $\G$ and the span of $\Pi$ densely generate $G$.

Whenever the metacurvature obstruction vanishes, the de~Rham complex of differential forms becomes a differential graded Poisson algebra. This is precisely the natural sufficient condition for the generalized Poisson brackets to descend to de~Rham cohomology. It is thus natural to ask what this gives in the cases we have been considering. The answer is disappointingly trivial.

\begin{prop}
If $\M$ is a compact Riemannian manifold with a compatible Poisson structure, then the induced Poisson bracket on de~Rham cohomology is $0$.
\end{prop}
\begin{proof}
Firstly, using the Leibniz identity \eqref{Leibniz}, we can extrapolate to an explicit formula for the generalized Poisson bracket of two differential forms $\sigma,\rho\in\Omega^*(\Mt)$,
\[
\{\sigma,\rho\} = \Pi^{AB} \Lie_{X_A}\sigma\wedge\Lie_{X_B}\rho
\mbox.
\]
If $\sigma$ and $\rho$ are closed, then this can be rewritten as
\[
\begin{split}
\{\sigma,\rho\} &= \Pi^{AB} d(X_A\inner\sigma)\wedge d(X_B\inner\rho) \\
&= d[(X_A\inner\sigma)\wedge \Lie_{X_B}\rho]
\mbox.
\end{split}
\]
The last expression in brackets is just a contraction of $\rho$ with $\D\sigma\in\Omega(\Mt,T\Mt)$. In this way, the formula descends to $\M$. If $\sigma,\rho\in\Omega^*(\M)$ are closed forms on $\M$, then $\D\sigma\in\Omega^*(\M,T\M)$ and so $\{\sigma,\rho\}$ is exact. 
\end{proof}

There is no obvious reason for the bracket on cohomology to vanish in greater generality; this question can be considered for any Poisson manifold that has a contravariant connection with vanishing torsion, curvature, and metacurvature. 
Unfortunately, the only other examples that I discuss in this paper are $S^2$ and $\R^3$; in those cases the bracket on cohomology vanishes simply because the cohomology is trivial.

\section{Spectral Triples}
I have shown that the compatibility conditions of Definition \ref{compatible} are necessary for the existence of a deformed noncommutative geometry which respects differential forms and integration. I have not shown if these conditions are sufficient. As I shall explain, this appears to be essentially true provided that the Poisson structure is itself suitably integrable.

I have tried to be as general as possible by not tying my arguments to a specific notion of noncommutative geometry more than necessary.  In order to discuss the sufficiency of my compatibility conditions, it is appropriate to be a bit more concrete. 

\subsection{Differential Forms}
\label{forms}
Connes \cite{con1} has given a very general recipe for constructing a differential graded algebra of ``noncommutative differential forms'' $\Omega^*_D(\A)$ from a spectral triple $(\A,\Hi,D)$. Given the Dirac operator and algebra of smooth functions on a compact spin manifold, this recovers the differential forms: 
\[
\Omega^*_D[\C^\infty(\M)] \cong \Omega^*(\M)
\mbox.
\]

The construction applies provided that the commutator of $D$ with any element of $\A$ is bounded. It begins by building a universal differential graded algebra from $\A$. Let $\Omega^0(\A)=\A$ and $\Omega^1(\A) \subset \A\otimes\A$ the kernel of the multiplication map $\A\otimes\A\to\A$. The differential begins with 
\[
d : \Omega^0(\A) \to \Omega^1(\A), \quad
a\mapsto 1\otimes a - a \otimes 1
\mbox.
\]
 Finally, $\Omega^k(\A)$ is the $k$-fold tensor product of $\Omega^1(\A)$ over $\A$, and $d$ is defined in general by the Leibniz identity. Note that $\Omega^k(\A) \subset \A^{\otimes (k+1)}$.

The formula $a_0\otimes a_1 \otimes \dots \otimes a_k \mapsto a_0Da_1D\dots D a_k$ defines a map from $\A^{\otimes(k+1)}$ to operators. Restricting this to $\Omega^k(\A)$ gives a representation by bounded operators $\rho : \Omega^*(\A) \to \Li(\Hi)$. The kernel of $\rho$ is neither a differential nor graded ideal. Instead define $J_0\subset \ker \rho \subset \Omega^*(\A)$ as the subspace spanned by homogeneous elements; then $J := J_0 + dJ_0$ is a differential graded ideal. Finally, $\Omega_D(\A) := \Omega^*(\A)/J$.

This doesn't require a spectral triple that could be reasonably regarded geometrically. For example, we could take $D=0$, in which case $\Omega^0_0(\A)=\A$ and $\Omega^{*>0}_0(\A)=0$.

\subsection{Axioms}
In \cite{con2}, Connes presented a system of axioms for a real spectral triple. This is the most completely and restrictively defined notion of noncommutative geometry. For reference, I summarize the axioms here. See \cite{gb-v-f} for the most detailed discussion.

An $n$-dimensional real spectral triple consists of: $\Hi$ a Hilbert space; $\A$ a $*$-algebra of bounded operators; $D$ an unbounded self-adjoint operator; $J$ an antiunitary operator; and $\gamma$ a $\Z_2$-grading operator. These satisfy the following axioms.

\emph{Dimension:} The resolvent $(D+i)^{-1}$ is a compact operator contained in the ideal $\Li^{n+}(\Hi)$.

\emph{Smoothness:} For any $a\in\A$, the commutator $[D,a]$ is bounded. Both $a$ and $[D,a]$ are in the domain of any power of the derivation $\delta$ defined by $\delta(a) := [\abs{D},a]$.

\emph{Reality:} For any $a,b\in\A$, $a^\op:= Ja^*J^{-1}$ commutes with $b$. $J^2=\pm1$, $JD=\pm DJ$, and $J\gamma= \pm \gamma J$, with the signs depending upon the dimension $n$ modulo $8$.

\emph{First order:} For any $a,b\in\A$, $a^\op$ commutes with $[D,b]$.

\emph{Finiteness:} The common domain $\Hi_\infty$ of all powers of $D$ is a finitely generated, projective $\A$-module. There is an $\A$-valued pre-Hilbert module inner product such that for $a\in\A$, $\psi,\chi\in\Hi_\infty$,
\[
\Tr_\omega \left( a\langle\psi\vert\chi\rangle_{\Hi_\infty} \abs{D}^{-n}\right)
= \langle \psi\rvert a^\op \lvert\chi\rangle_{\Hi}
\mbox.
\]

\emph{Orientation:} $\gamma$ is self-adjoint and commutes with $\A$. If $n$ is even then $\gamma D= - D\gamma$; if $n$ is odd then $\gamma=1$. There exists a Hochschild cycle $c\in Z_n(\A,\A\otimes \A^\op)$ such that $\rho(c)=\gamma$ for the representation $\rho: Z_*(\A,\A\otimes \A^\op) \to \Li(\Hi)$ defined by
\[
\rho(a_0\otimes b^\op \otimes a_1\otimes \dots \otimes a_n) 
= a_0 b^\op [D,a_1] \dots [D,a_n]
\mbox.
\]

\emph{\Poincare\ duality:} The Kasparov product with the $K$-homology class 
\[
[D]\in KR^n(\A\otimes \A^\op)
\]
is an isomorphism:
\[
\otimes_\A [D] : K_*(\A) \isom K^*(\A^\op) = K^*(\A)
\mbox.
\]

\subsection{Converse Construction}
\label{converse}
If a noncommutative deformation of the geometry of $\M$ exists, then we must have in particular a noncommutative deformation of the algebra of smooth functions on $\M$. Whether this exists is a more fundamental issue than compatibility with geometry, and in a way it is a separate question.  

The results of Section \ref{global} show that a compatible Poisson structure comes from a homogeneous Poisson manifold $G/\Gamma$.  This is the model that everything is constructed from. Assuming that $G/\Gamma$ can be suitably deformed, I will sketch the construction of real spectral triples for a deformation of $\M$.

I will focus on real spectral triples, so it is necessary to assume that $\M$ is a spin manifold. I assume that $\M$ is a compact, Riemannian spin manifold with a $\G$-covering $\Mt$, $\G\subset G$ is a discrete, cocompact subgroup, $\Pi\in\Wedge^2\g$ is an invariant Poisson bivector, and the action of $G$ on $\Mt$ extends to the spinor bundle. In order to satisfy the last condition, we might need to replace $G$ with some finite covering.

A geometric deformation of $\M$ exists if there is a suitable $G$-equivariant deformation of $G/\G$. I will not only assume that there exists a deformation $\BB$ of $\Ci(G/\G)$, but also:
\begin{enumerate}
\item
$\BB$ is a dense subalgebra of sections of a continuous field, $B$, of \cs-algebras over an interval.
\item
There is an action of $G$ on $\BB$ which extends the action of $G$ on $\Ci(G/\G)$, i.e., this is an equivariant deformation. As a deformation of $G$-modules, this is trivial.
\item
For every value of the parameter $\k$, let $\Bk\subset B_\k$ be the image of $\BB$ in $B_\k$. This is precisely the domain of the action of the universal enveloping algebra of $\g$ on $B_\k$.
\item
There exists a $G$-invariant tracial state $\tau_\k: B_\k\to\co$.
\end{enumerate}

It is elementary to construct a deformation of $G$ with respect to $\Pi$. In some cases, the techniques of Rieffel \cite{rie2} can then be used to construct a deformation of $G/\G$. However, it is not clear whether $G/\G$ can be suitably deformed in general.

\subsubsection{Algebra}
In preparation for constructing the deformed algebras, we can re-express $\M$ as
\[
\M \cong \Mt/\Gamma \cong \Mt \times_G G/\Gamma
\mbox.
\]
This means that the algebra of smooth functions $\Ci(\M)$ is naturally identified with the $G$-invariant $\Ci(G/\Gamma)$-valued smooth functions on $\Mt$:
\[
\Ci(\M) \cong \Ci(\Mt,\Ci[G/\Gamma])^G
\mbox.
\]
The (right) action of $g\in G$ is defined by the pullback by the right action of $g^{-1}$ on $\Mt$ and the the left action of $g$ on $G/\Gamma$.

We can mimic this construction in the noncommutative case substituting $\B_\k$ for $\Ci(G/\Gamma)$. There is still a right action of $G$ on $\B_\k$. The right action of $g\in G$ on $\Ci(\Mt,\B_\k)$ is given by simultaneously applying the action of $g$ on $\B_\k$ and pulling back by the action of $g^{-1}$ on $\Mt$.
With this action, we can define the algebra
\[
\Ak := \Ci(\Mt,\B_\k)^G
\mbox.
\]
This generalizes the ``twisting by a torus'' construction of Connes and Landi \cite{c-l} (see also \cite{c-dv1}).

\subsubsection{Hilbert Space}
Let $\Sn\to \M$ be the spinor bundle. The $\Ak$-module
\[
\Hi^\infty_\k:=\Gamma(\Mt,\Sn\otimes \B_\k)^G
\]
is an analogue of the smooth sections of $\Sn$ over $\M$, but in order to construct a spectral triple we need a Hilbert space analogous to the $L^2$-sections of $\Sn$. For this, we need an inner product.

The inner product of spinors comes fundamentally from a bundle homomorphism
\[
\bar\Sn\otimes\Sn \to \Wedge^n T^*\M
\mbox.
\]
For two spinor sections $\psi,\varphi\in\Gamma(\M,\Sn)$, the local inner product is $\bar\psi\varphi\in \Omega^n(\M)$. Integrating this over $\M$ gives the Hilbert inner product.

Now, for two $G$-invariant sections 
\[
\psi,\varphi \in \G(\Mt,\Sn\otimes \B_\k)^G
\]
the local inner product is $\bar\psi\varphi \in \Omega^n(\Mt,\B_\k)^G$. If we apply the $G$-invariant trace, this becomes
\[
\tau_\k(\bar\psi\varphi) \in \Omega^n(\Mt)^G
\mbox.
\]
It would be a mistake to try to integrate this over $\Mt$, since it does not fall off at all. Instead note that
\[
\Omega^n(\Mt)^G \subset \Omega^n(\Mt)^\Gamma \cong \Omega^n(\M)
\mbox.
\]
With this identification, we can integrate over $\M$ and define the Hilbert inner product as
\[
\langle\psi\mid\varphi\rangle := \int_\M \tau_\k(\bar\psi\varphi)
\mbox.
\]

It is worth observing how this gives the correct inner product at $\k=0$. In that case, $\bar\psi\varphi$ can already be identified with an $n$-form on $\M$. There is a singular foliation of $\M$ by the images of the $G$ orbits. The trace $\tau_0$ averages $\bar\psi\varphi$ over each of these leaves, which does not change the integral over $\M$.

\subsubsection{Dirac Operator}
The Dirac operator of $\M$ is defined on sections of $\Sn$. This extends trivially (and $G$-equivariantly) to sections of $\Sn\otimes \B_\k$ over $\Mt$. Restricting to $G$-invariant sections defines the Dirac operator $D_\k$ on $\Hi_\k$.

Identifying the algebras $\Ci(G/\G)$ and $\Bk$ as $G$-modules induces a unitary map from $L^2(\M,\Sn)$ to $\Hi_\k$ which intertwines the Dirac operators.
This is thus an isospectral deformation. We can identify the Hilbert spaces and regard the Dirac operator as constant. 
Because this is isospectral, the dimension axiom is trivially satisfied.
 
\subsubsection{Real Structure}
The classical real structure is given by an antilinear bundle automorphism $C:\Sn\to\Sn$. Combining this with the involution on $\B_\k$ gives the real structure $C\otimes*$ on $\Sn\otimes\B_\k$, and hence on $\Hi_\k$.

This leads to the obvious $\Ak$-bimodule structure on $\Hi_\k$. The various signs remain the same as in the commutative case. The first order axiom is easy to verify.

\subsubsection{Smoothness}
The action of $D_\k$ on $\G(\Mt,\Sn\otimes\Bk)^G$ can be rewritten partly in terms of the Lie algebra $\g$. The common domain of all powers of $D_\k$ is simply $\Hi^\infty_\k= \G(\Mt,\Sn\otimes\Bk)^G$. The finiteness axiom can be checked from this. Smoothness follows similarly.

\subsubsection{Differential Forms}
Let $\C^\infty_{\mathrm b}(\Mt,\Bk)$ denote the algebra of smooth $\Bk$-valued functions that are bounded in all derivatives. This is a tensor product of $\C^\infty_{\mathrm b}(\Mt)$ with $\Bk$. With $D\otimes 1$ and the Hilbert $B_\k$-module $L^2(\Mt,\Sn\otimes B_\k)$ this forms a spectral triple (in a slightly generalized sense). The algebra of noncommutative differential forms is simply 
\[
\Omega^*_{D\otimes 1}[\C^\infty_{\mathrm b}(\Mt,\Bk)] \cong 
\Omega^*_{\mathrm b}(\Mt,\Bk)
\mbox.
\]

Let $\rho_1: \Omega^*[\C^\infty_{\mathrm b}(\Mt)] \to \Li_{B_\k}[L^2(\Mt,\Sn\otimes B_\k)]$ be the representation used in the construction. The image consists of bounded-adjointable operators which act locally over $\Mt$. The same is true for the representation $\rho_2:\Omega^*(\Ak)\to \Li(\Hi_\k)$. Consequently, $\ker \rho_2=\ker \rho_1 \cap \Omega^*(\Ak)$. This means that $\Omega_{D_\k}(\Ak)$ is the differential-graded subalgebra of $\Omega^*_{\mathrm b}(\Mt,\Bk)$ generated by $\Ak$ in degree $0$. Thus 
\[
\Omega_{D_\k}(\Ak) \cong \Omega^*(\Mt,\Bk)^G
\mbox.
\]

With this, it is clear that the differential graded algebra $\Omega_{D_\k}(\Ak)$ is smoothly deformed from $\Omega^*(\M)$.

\subsubsection{Orientation}
By construction, the volume form $\tilde\epsilon\in\Omega^n(\Mt)$ is $G$-invariant. This means that it can be identified with $\tilde\epsilon\otimes 1 \in \Omega^n(\Mt,\Bk)^G \cong \Omega^n_{D_\k}(\Ak)$. Viewing this as an equivalence class in $\Omega^n(\Ak)$, some element of this class should be the Hochschild cycle with image $\gamma$ as required by the orientation axiom.

\subsubsection{\Poincare\ Duality}
This axiom is essentially impossible to check at this level of generality. It depends upon the stability of $K$-theory under the deformation from $\C^\infty(\M)$ to $\Ak$. This is not a general property, but there is a strong tendency for $K$-theory to be preserved in deformations, see \cite{ros}.

It is certainly plausible that \Poincare\ duality will be preserved in an isospectral deformation which preserves $K$-theory. For example, if the dimension is even, we can consider the intersection product of even $K$-theory classes which are determined by idempotents $e,e'\in \mathrm{Mat}_m(\Ak)$. Varying $\k$, the intersection product should vary continuously, but it is an integer and thus constant.

\section{Examples}
\label{examples}
\subsection{Torus}
The only compact example in two dimensions is a flat torus with a constant symplectic structure. In this case, $G$ is the 2-torus and $\G$ is trivial. 
The construction of the deformed geometry is the canonical example of noncommutative geometry.

More generally, any antisymmetric $n\times n$-matrix $\Pi$ defines a compatible Poisson structure on a flat $n$-dimensional torus, and there is a corresponding noncommutative torus deformation. The group $G$ is $\T^n$ in this case. Although $\Pi$ may be degenerate, the dimension of the group $T\normal G$ may still be larger than $\rk \Pi$. A simple example is $\R^3/\Z^3$ with $\pi = \partial_x\wedge(\partial_y+\sqrt{2}\partial_z)$. This is the effect I referred to in defining $T$ to be densely generated by the span of $\Pi$.  

\subsection{Flat Manifolds}
Let $\Mt = \T^4 = \R^4/\Z^4$ be the $4$-dimensional torus with coordinates $x^1,x^2,y^1,y^2$ (on $\R^4$) and symplectic form $\omega = dx^1\wedge dy^1 + dx^2\wedge dy^2$. The mapping 
\[
s : (x^1,y^1,x^2,y^2) \mapsto (y^1,-x^1,x^2+\tfrac14,y^2)
\]
preserves $\omega$ and generates a free action of $\Z_4$ on $\Mt$. Define $\M:=\Mt/\Z_4$. 

The abelian group of Poisson isometries generated by the span of $\Pi$ is the full group of translations, $\T^4$. The covering group is $\G\cong\Z_4$ generated by $s$. The group generated by these is the semidirect product $G= \Z_4 \ltimes \T^4$; the generator $s\in\Z_4$ acts on $\T^4$ by a quarter rotation in the $x^1$-$y^1$-plane. In this case, $\M$ is of the locally homogeneous form $\M \cong H\backslash G /\G$ with $H\cong \G \cong \Z_4$.

As always, we should construct a noncommutative deformation by first deforming the homogeneous model $G/\G\cong \T^4$. This gives a noncommutative 4-torus. The algebra $\Bk$ is generated by 4 unitaries $U_1$, $U_2$, $V_1$, and $V_2$ with the relations
\[
\begin{split}
U_1V_1 &= e^{(2\pi)^2i\k}V_1U_1, & U_2V_2, &= e^{(2\pi)^2i\k}V_2U_2, 
\end{split}
\]
and all other pairs of generators commuting. Corresponding to $s$ is the automorphism $s^*$ defined by: 
$s^*(U_1) = V_1$, $s^*(V_1) = U_1^{-1}$, $s^*(U_2) = iU_2$, and $s^*(V_2)=V_2$. $\Ak$ is the $s^*$-invariant subalgebra of $\Bk$.

This example is reminiscent of the classification of flat manifolds: Any flat compact Riemannian manifold is of the form $\M=\T^n/\G$ with $\G$ a finite group of free isometries of the flat torus (see \cite{cha}). If there exists a constant $\G$-invariant Poisson structure on $\T^n$, then this descends to a Poisson structure on $\M$. Choosing $\Mt=\T^n$, we can take $G=\G\ltimes \T^n$ and this is another locally homogeneous example $\M \cong \G\backslash G /\G$.
The standard quantization of $\T^n$ is equivariant under all linear transformations preserving the Poisson structure, therefore this gives a deformation of $\M$.

In general, if $\M$ is a compact Riemannian manifold with a compatible symplectic structure, then it is flat and it must be of this form.

\subsection{Heisenberg Manifolds}
The Heisenberg group $H$ is the group of real $3\times3$ matrices of the form
\[
\begin{pmatrix}
1 & x & z\\
0 & 1 & y \\
0 & 0 & 1
\end{pmatrix}
\mbox.
\]
Let $\Lambda\subset H$ be the subgroup of matrices with integer entries. The 3-dimensional nil-manifold (or ``Heisenberg manifold'') is $\Nil := H/\Lambda$.

Let $\theta$ be an irrational number. The Poisson bivector 
\[
\Pi = (\partial_x+\theta\partial_y)\wedge\partial_z
\]
is bi-invariant on $H$. In particular it is $\Lambda$-invariant and defines a Poisson structure on $\Nil$. This is a regular Poisson structure. The leaves of the symplectic foliation are dense. $\Nil$ is the total space of a circle bundle over the torus $\T^2$; the symplectic foliation is the inverse image of a Kronecker foliation of the torus.

Obviously $\Nil=H/\Lambda$ describes $\Nil$ in the desired form, but this is not the minimal description as in Thm.~\ref{structure}. The universal covering $H$ is not the minimal covering of $\Nil$ with the Poisson bivector given by global Killing vectors. Instead, we can take $\Mt = H/\Z$ where $\Z$ is embedded in $H$ as the subgroup of matrices
\[
\begin{pmatrix}
1 & 0 & z\\
0 & 1 & 0\\
0 & 0 & 1
\end{pmatrix}
\]
with $z\in\Z$.

The covering group is $\G = \Lambda/\Z \cong \Z^2$. The image of $\Pi$ in the Heisenberg Lie algebra is 2-dimensional, spanned by
\[
\begin{pmatrix}
0&1&0\\0&0&\theta\\0&0&0
\end{pmatrix}
\text{ and }
\begin{pmatrix}
0&0&1\\0&0&0\\0&0&0
\end{pmatrix}
\mbox.
\]
This generates an abelian subgroup of Poisson isometries $T\cong \R\times\T^1$.

This covering manifold $\Mt$ is itself a group $G=H/\Z$ because $\Z\normal H$. Both $\Gamma$ and $T$ are contained in the group of right translations. In fact $T\,\G\subset G$ is dense. This is easiest to see by looking at the quotient $G/\T^1 \cong \R^2$. $\G$ maps injectively to the integer lattice $\Z^2\subset\R^2$. $T$ maps to a 1-dimensional subgroup, a line of irrational slope. Together, these densely generate $\R^2$.

As a Poisson manifold, $\Nil$ is homogeneous. The left action of $G$ on $\Nil\cong G/\Z^2$ preserves the Poisson structure. However, it is not homogeneous as a Riemannian manifold. A compatible Riemannian metric is given by any right-invariant metric on $G$. This cannot be left-invariant as well. 

Rieffel constructed an equivariant deformation quantization of $\Nil$ in \cite{rie1}. In fact, this was one of the very first examples of strict deformation quantization. The other Heisenberg manifolds discussed there are simply finite quotients of $\Nil$. They can be seen as locally homogeneous examples subordinate to $\Nil$.
Chakraborty and Sinha \cite{c-s1} have constructed and analyzed spectral triples for these examples.

\section{Conclusions}
With the assumptions I have made about what constitutes noncommutative geometry, I have shown that noncommutative deformations are remarkably restricted.

The reader might wonder if the analysis here was all necessary. After all, the conclusion of Thm.~\ref{local.structure} is quite simple. The reason for the complexity of the proof is that compactness leads to great simplification, but in a very indirect way. The noncompact solutions of the compatibility conditions are more complicated. For instance, consider $3$-dimensional Euclidean space with Cartesian coordinates $x,y,z$, and a Poisson bracket defined by $\{x,y\}=1$, $\{x,z\}=y$, and $\{y,z\}=-x$. This satisfies the compatibility conditions, but the Poisson bivector cannot be decomposed into products of Killing vectors. This example can be noncommutatively deformed, which suggests that there does not exist another local obstruction that would have simplified the analysis in Section \ref{local}.

Another criticism is that my compatibility conditions are too restrictive. Although I think that my assumptions are well founded, some variation is possible. Different assumptions may lead to weaker conclusions, but I think that this has been a suitable starting point. The techniques I have used here may be useful for analyzing other scenarios.

The \Dabrowski-Sitarz example of a spectral triple for the \Podles\ sphere is an example of noncommutative geometry where the orientation axiom is not satisfied. This is partly because the homological dimension of the \Podles\ sphere is $0$ rather than $2$. This sudden dimension drop is related to the fact that integration cannot be smoothly deformed to a trace. This example strongly suggests that it would be interesting to discard my compatibility condition $0=d(\pi\inner\epsilon)$ between the Poisson structure and volume form. Unfortunately, this condition played a key role in the simplifications leading to the main results here. Nevertheless, the technique of using symplectic realizations to tame the metacurvature still applies. There may be other tricks that would make this problem tractable.

Some constructions of noncommutative differential forms for the fuzzy sphere and quantum groups have avoided these obstructions by using differential forms that do not correspond to classical differential forms. Instead of deforming $1$-forms, they deform $1$-forms plus some other ``junk''. This approach seems unpleasantly \emph{ad hoc} to me, but it is probably still possible to analyze this situation with my techniques.

The construction of spectral triples on noncommutatively deformed spaces is an active area of research. Theses examples are interesting, but in my view they cannot be considered as deformed geometries. For instance, the spectral triple for $\SU_q(2)$ constructed by Chakraborty and Pal \cite{c-p1} treats $\SU_q(2)$ as a quantum group in its own right, rather than as a deformation of $\SU(2)$.

Physically, these results rule out the idea of noncommutatively deformed 4-di\-men\-sional space-time. At the same time, they spell out a possible structure for noncommutative extra dimensions.

I have essentially classified noncommutative deformations of compact Riemannian manifolds in terms of the structure $(G,\G,\Pi)$ where $G$ is a Lie group, $\G\subset G$ is discrete and cocompact, $\Pi\in\Wedge^2\g$ is an invariant Poisson structure, and $G$ is densely generated by $\G$ and the Lie algebra ideal spanned by $\Pi$. What is missing is a better understanding of this structure. For example, the Lie algebra $\g$ is certainly not semisimple (it has a nontrivial ideal) but in all the examples that I know of, $\g$ is actually nilpotent; it is not apparent whether $\g$ is always nilpotent.

\end{document}